\documentclass[11pt]{amsart}
\usepackage{graphicx}
\usepackage{amsmath,amssymb,mathrsfs,extpfeil,latexsym}
\usepackage{algorithm,algpseudocode}
\usepackage{indentfirst}
\usepackage{enumerate}
\usepackage{float}
\usepackage{epsfig,subcaption}
\usepackage{colortbl}
\usepackage{hyperref}
\usepackage[title]{appendix}

\begin{document}

\title{Model predictive control with random batch methods for a guiding problem}

%

\author[D. Ko]{Dongnam Ko}
\address{
\begin{flushleft}
Dongnam Ko
\end{flushleft}
\begin{flushleft}
Chair of Computational Mathematics, Fundaci\'on Deusto, University of Deusto, 48007 Bilbao, Basque Country, Spain
\end{flushleft}
\begin{flushleft}
dongnamko@deusto.es
\end{flushleft}
}

\author[E. Zuazua]{Enrique Zuazua}
\address{
\begin{flushleft}
Enrique Zuazua
\end{flushleft}
\begin{flushleft}
Chair in Applied Analysis, Alexander von Humboldt-Professorship, Department of Mathematics, Friedrich-Alexander-Universit\"at Erlangen-N\"urnberg,
91058 Erlangen, Germany,
\end{flushleft}
\begin{flushleft}
and
\end{flushleft}
\begin{flushleft}
Chair of Computational Mathematics, Fundaci\'on Deusto, University of Deusto, 48007 Bilbao, Basque Country, Spain,
\end{flushleft}
\begin{flushleft}
and
\end{flushleft}
\begin{flushleft}
Departamento de Matem\'aticas, Universidad Aut\'onoma de Madrid, 28049 Madrid, Spain 
\end{flushleft}
\begin{flushleft}
enrique.zuazua@fau.de
\end{flushleft}
}

\newtheorem{theorem}{Theorem}[section]
\newtheorem{lemma}{Lemma}[section]
\newtheorem{corollary}{Corollary}[section]
\newtheorem{proposition}{Proposition}[section]
\newtheorem{example}{Example}[section]
\newtheorem{remark}{Remark}[section]
\newtheorem{definition}{Definition}[section]


\thanks{\textbf{Acknowledgment.} This project has received funding from the European Research Council (ERC) under the European Union's Horizon 2020 research and innovation programme (grant agreement No. 694126-DyCon). 
The work of the second author has been funded by the Alexander von Humboldt-Professorship program, the European Union's Horizon 2020 research and innovation programme under the Marie Sklodowska-Curie grant agreement No.765579-ConFlex, grant MTM2017-92996-C2-1-R COSNET of MINECO (Spain), ELKARTEK project KK-2018/00083 ROAD2DC of the Basque Government, ICON of the French ANR and Nonlocal PDEs: Analysis, Control and Beyond, AFOSR Grant FA9550-18-1-0242, and Transregio 154 Project *Mathematical Modelling, Simulation and Optimization using the Example of Gas Networks* of the German DFG}

\begin{abstract}
We model, simulate and control the guiding problem for a herd of evaders under the action of repulsive drivers. The problem is formulated in an optimal control framework, where the drivers (controls) aim to guide the evaders (states) to a desired region of the Euclidean space.
The numerical simulation of such models quickly becomes unfeasible for a large number of interacting agents.
To reduce the computational cost, we use the Random Batch Method (RBM), which provides a computationally feasible approximation of the dynamics. 
At each time step, the RBM randomly divides the set of particles into small subsets (batches), considering only the interactions inside each batch.
Due to the averaging effect, the RBM approximation converges to the exact dynamics as the time discretization gets finer. 
We propose an algorithm that leads to the optimal control of a fixed RBM approximated trajectory using a classical gradient descent. 
The resulting control is not optimal for the original complete system, but rather for the reduced RBM model. 
We then adopt a Model Predictive Control (MPC) strategy to handle the error in the dynamics.
While the system evolves in time, the MPC strategy consists in periodically updating the state and computing the optimal control over a long-time horizon, which is implemented recursively in a shorter time-horizon.
This leads to a semi-feedback control strategy. Through numerical experiments we show that the combination of RBM and MPC leads to a significant reduction of the computational cost, preserving the capacity of controlling the overall dynamics.
\end{abstract}

\maketitle

\section{Introduction}
\setcounter{equation}{0}

Control problems on collective behavior systems have received huge attention recently, \cite{bailo_optimal_2018,bongini_sparse_2015,fax2004information,ma_finite-time_2017}, due to the growing needs in applications.
In the context of interacting particle systems, a collective behavior system  can be viewed as a decentralized control system, \cite{carrillo_analytical_2018,dorfler_synchronization_2012,motsch_new_2011}, where each individual reacts to its interacting agents in order to achieve its own nature or desire. In such self-organized or self-driven systems, a centralized control can play as an external interference boosting the collective behavior to evoke the desired dynamics

 As indicated in literature, \cite{bongini_mean-field_2017}, from a practical viewpoint, the control of a complex system is commonly designed to manipulate a small portion of the particles,  \cite{caponigro_sparse_2013,piccoli_sparse_2019}, or preselected agents, playing the role of leaders  or informed agents, \cite{fax2004information,porfiri2008criteria}.
 
One of the most relevant examples is the shepherding problem, where a few shepherd dogs are required to steer a herd of sheep.
Many attempts have been made in various contexts to understand how the dogs can affect a group of sheep, \cite{strombom_solving_2014}, and steer them to the desired region, \cite{jyh-ming_lien_shepherding_2004}. 

Our interest in this paper consists in modeling, deriving  and  simulating optimal control strategies for guiding problems, where a small number of repelling agents (drivers), which play the role of controls, have to guide a herd of flocking agents (evaders), that we interpret as states,  toward a given desired area. As a basis, we follow the formulation based on the   \emph{guidance-by-repulsion} paradigm, \cite{escobedo_optimal_2016,ko2019GBR}. 

The problem with one driver and one evader has been addressed in \cite{escobedo_optimal_2016}. Its long-time behavior and controllability properties were analyzed in \cite{ko2019GBR}. In that simplified setting numerical simulations show that the guiding problem is very sensitive to the motion of drivers, the optimization process being highly non-convex.

The main novelty of this article is that we consider an arbitrary number of evaders and guiders, developing new control and computational techniques allowing to handle large systems rather precisely. Of course, with practical applications in mind, the number of evaders is large, but one aims at controlling the herd efficiently by a reduced number of drivers.

The first difficulty encountered when addressing these issues is the computational complexity and cost of the forward dynamics, which dangerously increases when implementing the iterative algorithms needed to compute optimal control strategies. In particular, proceeding as in \cite{ko2019GBR}, the optimal control with four drivers and sixteen evaders in a time-horizon $T=40$, with time step $\Delta t = 0.01$, requires nearly two hours in a laptop computer with CPU i5-4258U 2.4GHz and RAM DDR3L 8GB 1600MHz, operated in Matlab.

In this paper, to get a reliable control strategy with an affordable computational cost, we suggest an algorithm (see Algorithm \ref{alg:MPC-RBM} in Section \ref{sec:pre}), based on a classical optimal control formulation, but combining critically two main ingredients: First, the Random Batch Method (RBM), \cite{jin_random_2020}, which provides and approximation of the dynamics of the interacting particle system at a small computational cost and, second, Model Predictive Control (MPC), [18], to stabilize the approximation error, \cite{garcia_model_1989}.

The RBM is an approximation method for large systems of particles, particularly suitable in the context of collective dynamics when the individuals are not distinguishable. If a system has binary interactions between $N$ particles, the number of interactions is in the order of $O(N^2)$, which makes the computational cost of the forward dynamics extremely expensive.
To tackle this difficulty, instead of computing the whole interactions, for a given $1<P<N$, the RBM approximates dynamics out of $O(NP)$  interactions. 
More precisely, for a small duration of time, we split the set of particles into random small subsets (batches) which contain, at most, $P$ particles. Then, one only considers the interactions within each batch, ignoring the interaction between batches. In the next time interval, to average the random effect in time, we again choose batches independently.

Therefore, from the all-to-all interacting particle system, the RBM produces a deterministic networked model which periodically switches the network structure. Thanks to the random choices the reduced RBM model approximates the original time evolution properly. The approximation error of the RBM is analyzed in \cite{jin_random_2020}, and briefly presented in Section \ref{sec:2-6}. This is based on the Law of Large Numbers, which guarantees the convergence to the original dynamics as $\Delta t$ decreases to zero. The approximation error is also tested with numerical simulations in Section 3.1, which draw the $95\%$ confidence intervals from 200 independent simulations with $T = 10$ and $\Delta t = 0.01$. 

In this paper, we solve the optimal control problem on the simplified RBM model using standard gradient descent methods. But, of course, this does not lead to an accurate control of the complete dynamics.
The performance of the control depends on the approximation error, which accumulates in time. 
One can overcome this difficulty adopting the viewpoint of Model Predictive Control (MPC),  \cite{garcia_model_1989,grune2017nonlinear,nikolaou2001model}.

MPC allows to manage two important difficulties arising in complex systems; first,  the exact input-output equations  may be unknown or only partially known due to, for instance, noisy phenomena  making it difficult or even impossible  to predict the exact dynamics. Second, the system may have a sensitive behavior over a long time that a simple feedback controller may not regulate. For example, a greedy control strategy may  push the evaders toward the target too strongly, and this may result on making evaders pass over the target in the long run.

The basic idea of MPC is to build and approximated reduced model, the RBM one in the present context,   and to compute its open-loop optimal control within a long-time horizon. While the system evolves in time, the controls designed that way may fail to steer the complete system efficiently, since this open-loop control has no means of detecting errors or unexpected perturbations on the dynamics and to cope with the intrinsic gap between the complete and the reduced dynamics. Hence, after a short time, the state of the system has to be rechecked to recompute the optimal control with the updated initial data. In this iterative way, one can force  the control  to adapt to the true dynamics of the system.
A more detailed description of the RBM and MPC is presented in Section \ref{sec:pre}.

As shown in this paper, the adequate combination of the RBM and MPC may lead to a significant reduction of the computational cost preserving the efficiency of the control strategy to steer the complete dynamics. For instance, in one of the examples we describe in Section \ref{sec:simul}, our method shows an $84\%$ reduction of the computational time on the problem of $36$ evaders and $2$ drivers with $T = 4$ and $\Delta t = 0.01$, while the running cost (performance of control) only differs by about $0.3\%$, compared to the standard optimal control problem.

The rest of this paper is organized as follows. In Section \ref{sec:2-1}, we formulate the guiding problem as an optimal control one. From Section \ref{sec:2-2} to Section \ref{sec:2-4}, we present the preliminaries on the RBM and MPC. Then, the detailed procedure of how to combine MPC-RBM to build our algorithm is described in Section \ref{sec:2-5}. The error estimation on this algorithm is an open problem, as described in Section \ref{sec:2-6}. The simulations are presented in Section \ref{sec:simul} to test computational costs and approximation errors.
 Finally, in Section \ref{sec:final}, we discuss our results and present some final remarks and open problems arising in this field.

\section{The MPC-RBM algorithm}\label{sec:pre}

In this section, we present in detail the MPC algorithm using the RBM. We first present the guiding optimal control problem. Then the combined MPC-RBM algorithm is described.

\subsection{Optimal control formulation on the guiding problem}\label{sec:2-1}

Let ${\mathbf x} = ({\mathbf x}_1,\ldots,{\mathbf x}_N) \in \mathbb R^{Nd}$ and ${\mathbf v} = \dot{{\mathbf x}} \in \mathbb R^{Nd}$ be the positions and velocities of $N$  Newtonian particles moving in a $d$-dimensional space, namely, evaders. The control is indirectly introduced through $M$  particles represented by ${\mathbf y} = ({\mathbf y}_1,\ldots,{\mathbf y}_M) \in \mathbb R^{Md}$, called drivers. We assume that the evaders interact with other evaders and the drivers as in the following collective behavior system: for $t \geq 0$,
\begin{equation}\label{eq:GBR}
\begin{aligned}
\begin{cases}
\dot{\mathbf x}_i = {\mathbf v}_i,\quad i=1,\ldots,N,\\
\dot{\mathbf v}_i = \displaystyle\frac{1}{N-1}\sum_{k=1,k\neq i}^N a({\mathbf x}_k-{\mathbf x}_i)({\mathbf v}_k-{\mathbf v}_i)
+ \frac{1}{N-1}\sum_{k=1,k\neq i}^N g({\mathbf x}_k-{\mathbf x}_i)({\mathbf x}_k-{\mathbf x}_i)\\
\hspace{2em}
\displaystyle- \frac{1}{M}\sum_{j=1}^M f({\mathbf y}_j-{\mathbf x}_i)({\mathbf y}_j-{\mathbf x}_i),~ i=1,\ldots,N,\\
\dot {\mathbf y}_j = {\mathbf u}_j(t), \quad j=1,\ldots,M\\
{\mathbf x}_i(0) = {\mathbf x}_i^0,\quad {\mathbf v}_i(0) = {\mathbf v}_i^0, \quad {\mathbf y}_j(0) = {\mathbf y}_j^0.
\end{cases}
\end{aligned}
\end{equation}
We assume that the controls ${\mathbf u}(t) = ({\mathbf u}_1(t),$ $\ldots,$ ${\mathbf u}_M(t)) \in \mathbb R^{Md}$ completely determine the dynamics of the drivers, as a means to indirectly influence the dynamics of the evaders. 

The nonlinearities $a(\cdot)$, $f(\cdot)$ and $g(\cdot)$ entering in the dynamics are assumed to be smooth and positive, except for $g(\cdot)$ that may have negative values as in the context of intermolecular forces, to avoid collisions between evaders \cite{carrillo_sharp_2017,cucker_avoiding_2010}. As a concrete example and for the purpose of developing the numerical experiments, we define $a(\cdot)$, $f(\cdot)$ and $g(\cdot)$ as follows: for $i,k = 1,\ldots,N$ and $j=1,\ldots, M$,
\begin{equation}\label{eq:agf}
\begin{aligned}
a({\mathbf x}) &:= 1,\quad
f({\mathbf x}) := 4\exp(-8|{\mathbf x}|^2)\quad \text{and}\quad
g({\mathbf x}) := \begin{cases} \displaystyle 2\left( 1 -  \frac{1}{3\sqrt{N}|{\mathbf x}|^2} \right) \quad \text{if} \quad {\mathbf x} \neq 0,\\
0 \hspace{8.5em} \text{otherwise}.\end{cases}
\end{aligned}
\end{equation}

Note that \eqref{eq:agf} contains an unbounded function $g(\cdot)$, but  local existence and uniqueness is guaranteed for initial data satisfying ${\mathbf x}_i^0 \neq {\mathbf x}_j^0$ for any $i\neq j$. The dynamical properties of the system \eqref{eq:GBR} are discussed in \cite{gade_herding_2015} with a similar formulation of interactions. The global existence of solutions for \eqref{eq:agf} follows similar arguments to \cite{cucker_avoiding_2010}.

The qualitative properties of the nonlinearities in \eqref{eq:agf} are chosen to reflect the following main features of the ``drivers-evaders" interactions (see also \cite{gade_herding_2015,ko2019GBR,pinnau_interacting_2018,strombom_solving_2014} for specific guiding problems):
\begin{itemize}
\item
The evaders are influenced by a repulsive force $f(\cdot)$ from each driver,  its strength being decreasing as the distance increases.
\item
Each evader has positional interactions with other evaders given by the function $g(\cdot)$, which reflects their aim to remain close together. To prevent the collisions between evaders, we also set $g(\cdot)$ to have a strong negative value when a pair of evaders is too close. Overall, all the evaders have the tendency to gather in a disk with a diameter $\sim 0.5$.
\item
The evaders also interact each other, through the term $a(\cdot)$, to align the velocities to a common value. This plays the role of friction to the mean velocity of evaders, which reduces the oscillatory behavior arising from $g(\cdot)$.
\end{itemize}

Compared to  herding, \cite{jyh-ming_lien_shepherding_2004,strombom_solving_2014}, and  flocking problems, \cite{cucker_emergent_2007,ha2009simple}, the above system \eqref{eq:GBR} contains both the velocity alignment $a(\cdot)$ and the positional potential interactions $g(\cdot)$. The combination of positional and velocities' interactions is suggested in \cite{gade_herding_2015,park_cucker-smale_2010,tanner_flocking_2007}, which is used to model the dynamics of birds.

Our control objective is to guide the evaders to a desired region by the locomotion induced by the drivers in a given time horizon. In order to formulate this problem in the context of optimal control, we define the cost function as follows:
\begin{equation}\label{eq:cost}
\begin{aligned}
J({\mathbf u}) := \int_0^T &\left[ \frac{\alpha_1}{N} \sum_{k=1}^N |{\mathbf x}_k-{\mathbf x}_f|^2
 + \frac{\alpha_2}{M} \sum_{j=1}^M |{\mathbf u}_j|^2 + \frac{\alpha_3}{M} \sum_{j=1}^M |{\mathbf y}_j-{\mathbf x}_f|^2\right] dt.
\end{aligned}
\end{equation}
This functional takes account, in particular,  of the running cost of the distances to the target point ${\mathbf x}_f \in \mathbb R^d$. The positive constants $\alpha_1,\alpha_2,\alpha_3$ allow to regulate the weight of each of the terms entering in the cost. In practice, $\alpha_1$ is taken to be large compared to $\alpha_2$ and $\alpha_3$. 
Involving the running cost of the control by means of $\alpha_2$ is rather standard, to avoid unfeasibly large controls. On the other hand, the term involving $\alpha_3$ prevents the drivers not to go very far from the evaders. The target position ${\mathbf x}_f$ is chosen as a reference point since we expect that, tracking the evaders, the drivers will also spend most of the time near ${\mathbf x}_f$.

We expect that, minimizing this functional, the optimal control of the system will guide the evaders toward the region near ${\mathbf x}_f$ and capture them for a long time. Of course, the efficiency on doing so will depend on the length of the time-horizon and the value of the parameters $\alpha_j$, $j=1, 2, 3$.

The particular case of one driver and one evader was analyzed in \cite{escobedo_optimal_2016,ko2019GBR}. In that case, it was proved that  there exist a final time $T$ and a control function ${\mathbf u}_1 \in L^\infty((0,T),\mathbb R^{d})$ satisfying ${\mathbf x}_1(T) = {\mathbf x}_f$. But, of course, achieving such controllability results is much harder for multiple drivers and evaders. This is a very interesting open analytical problem, out of the scope of the present article.

To analyze a large number of particles (evaders), one of the well-known methods is the mean-field limit \cite{golse_mean-field_2003,ha2009simple}. It considers the distribution of particles instead of the whole trajectories, in the form of kinetic (or transport) equations. Hence, the evaders need to be undistinguishable; the interactions in \eqref{eq:GBR} only depend on the relative positions and velocities, not the index $i$ or $k$. In addition, if the interactions are bounded and smooth, then the mean-field limit strategy can be applied, \cite{bongini_mean-field_2017,burger_controlling_2016,pinnau_interacting_2018}.

In this paper, we study the particle description of the guiding problem due to its simplicity and generality. 
Our method, combining  MPC and RBM, can be easily generalized to the mean-field control problems and other formulations since the corresponding dynamics can be approximated by the particle model through the characteristic equations or the space discretization.

\subsection{The RBM approximation for the forward dynamics}\label{sec:2-2} 

We use the RBM in \cite{jin_random_2020} to reduce the computational cost to simulate the time evolution of the interacting particle system \eqref{eq:GBR}. 

We proceed as follows. 
\begin{itemize}
\item The control time interval $[0,T]$ is fixed. We fix a short duration of time $\Delta t>0$ and the discrete times $t_n := n \Delta t$. 

\item
Then, for each interval $[t_n,t_{n+1}]$, we independently choose $N/P$ random batches ${\mathcal C}_n^1, {\mathcal C}_n^2, \ldots, {\mathcal C}_n^{N/P}$ with size $P$ as a partition of the index set $\{1,\ldots,N\}$. 

\item
For a given integer $P>1$, that we assume is divisor of  $N$, in each time sub-interval  $[t_n,t_{n+1}]$, we divide the set of $N$ evaders into   $N/P$ batches of size $P$  that are chosen randomly. In case $N/P$ is not an integer, the last batch will have size less than $P$. 

\item In this way, each evader belongs to one and only one batch in each sub-interval $[t_n,t_{n+1}]$. We build a reduced dynamics so that we only compute the interactions inside each batch. Then, the number of interactions between the evaders decreases to $O(NP)$ from the all-to-all number of interactions $O(N^2)$.

\end{itemize} 
 
The resulting dynamics is of switching nature. In each time sub-interval,  the RBM rearranges the distribution of evaders into batches and modifies the dynamics.

We denote by ${\mathcal C^p_n}$, $p=1,..., N/P$ the $N/P$ batches corresponding to each subinterval $[t_n,t_{n+1}]$. 

For each index $i$ identifying an evader, by ${\mathcal C^{p(i)}_n}$ we denote the batch that, in the time interval $[t_n,t_{n+1}]$, contains this evader. Then, the approximated dynamics on the $i$th evader ($ i \in {\mathcal C^{p(i)}_n}$) can be formulated as follows: for $t \in [t_n,t_{n+1}]$,
\begin{equation}\label{eq:sGBR} 
\begin{cases}
\dot{\mathbf x}_i = {\mathbf v}_i, \qquad  \qquad t \in [t_n,t_{n+1}], \quad i \in {\mathcal C^{p(i)}_n},\\
\dot{\mathbf v}_i = \displaystyle\frac{1}{P-1}\sum_{k \in {\mathcal C^{p(i)}_n}\setminus\{i\}} a({\mathbf x}_k-{\mathbf x}_i)({\mathbf v}_k-{\mathbf v}_i)
 + \frac{1}{P-1}\sum_{k \in {\mathcal C^{p(i)}_n}\setminus\{i\}} g({\mathbf x}_k-{\mathbf x}_i)({\mathbf x}_k-{\mathbf x}_i)\\
\hspace{2em}\displaystyle - \frac{1}{M}\sum_{j=1}^M f({\mathbf y}_j-{\mathbf x}_i)({\mathbf y}_j-{\mathbf x}_i),
\qquad t \in [t_n,t_{n+1}], \quad i \in {\mathcal C^{p(i)}_n},\\
\dot {\mathbf y}_j = {\mathbf u}_j(t), \qquad j=1,\ldots,M,\\
{\mathbf x}_i(t_n) = {\mathbf x}_i^n,\quad {\mathbf v}_i(t_n) = {\mathbf v}_i^n, \qquad i \in {\mathcal C^{p(i)}_n},\\
 {\mathbf y}_j(t_n) = {\mathbf y}_j^n,  \qquad j=1,\ldots,M,\
\end{cases}
\end{equation}
In \eqref{eq:sGBR}, the scaling constant $1/(P-1)$ is taken to average the interactions instead of $1/(N-1)$, since there are $P-1$ interacting evaders for the $i$th evader in ${\mathcal C^{p(i)}_n}$. It averages the interactions on each evader as in the original system \eqref{eq:GBR}.

Note that the system \eqref{eq:sGBR} has a switching nature since it  is defined differently on each interval $[t_n,t_{n+1}]$. 
 From the initial data ${\mathbf x}_i^0$, ${\mathbf v}_i^0$ and ${\mathbf y}_j^0$, the RBM computes \eqref{eq:sGBR} with the batches ${\mathcal C}_0^1, {\mathcal C}_0^2, \ldots, {\mathcal C}_0^{N/P}$ over $t \in [0,t_1]=[t_0,t_1]$. Then, the final data ${\mathbf x}_i(t_1)$, ${\mathbf v}_i(t_1)$ and ${\mathbf y}_j(t_1)$ are taken as  the initial data for the next time interval, ${\mathbf x}_i^1$, ${\mathbf v}_i^1$ and ${\mathbf y}_j^1$ and one switches to the new section of batches ${\mathcal C}_1^1, {\mathcal C}_1^2, \ldots, {\mathcal C}_1^{N/P}$. Following the iterative calculations from $n=0$ to $N_T := \left \lfloor{T/\Delta t}\right \rfloor$, we can get the approximated time evolution in $[0,T]$. 
 
 \begin{figure*}[ht]
  \centering
  {
    \includegraphics[width=0.32\textwidth]{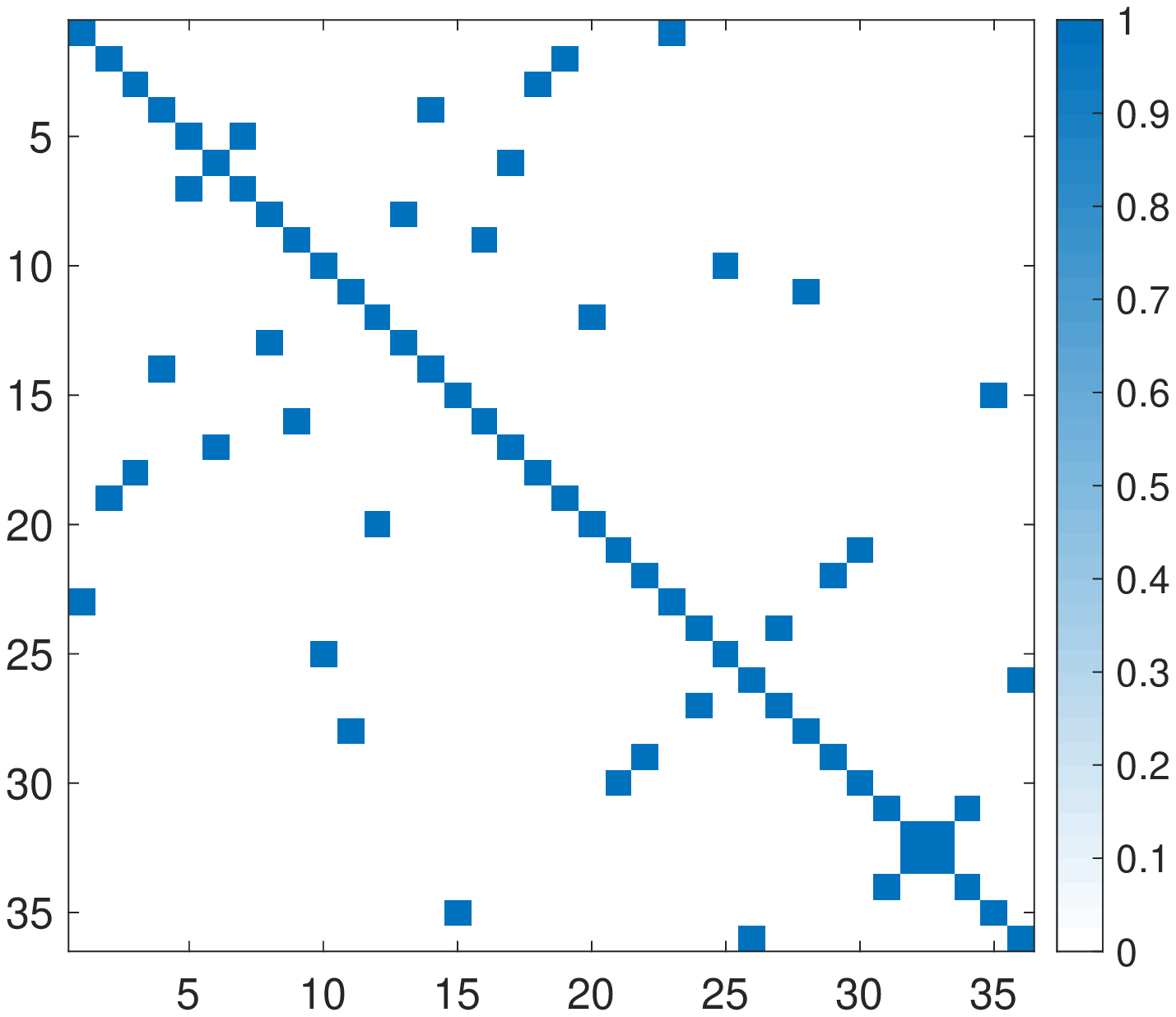}
    \includegraphics[width=0.32\textwidth]{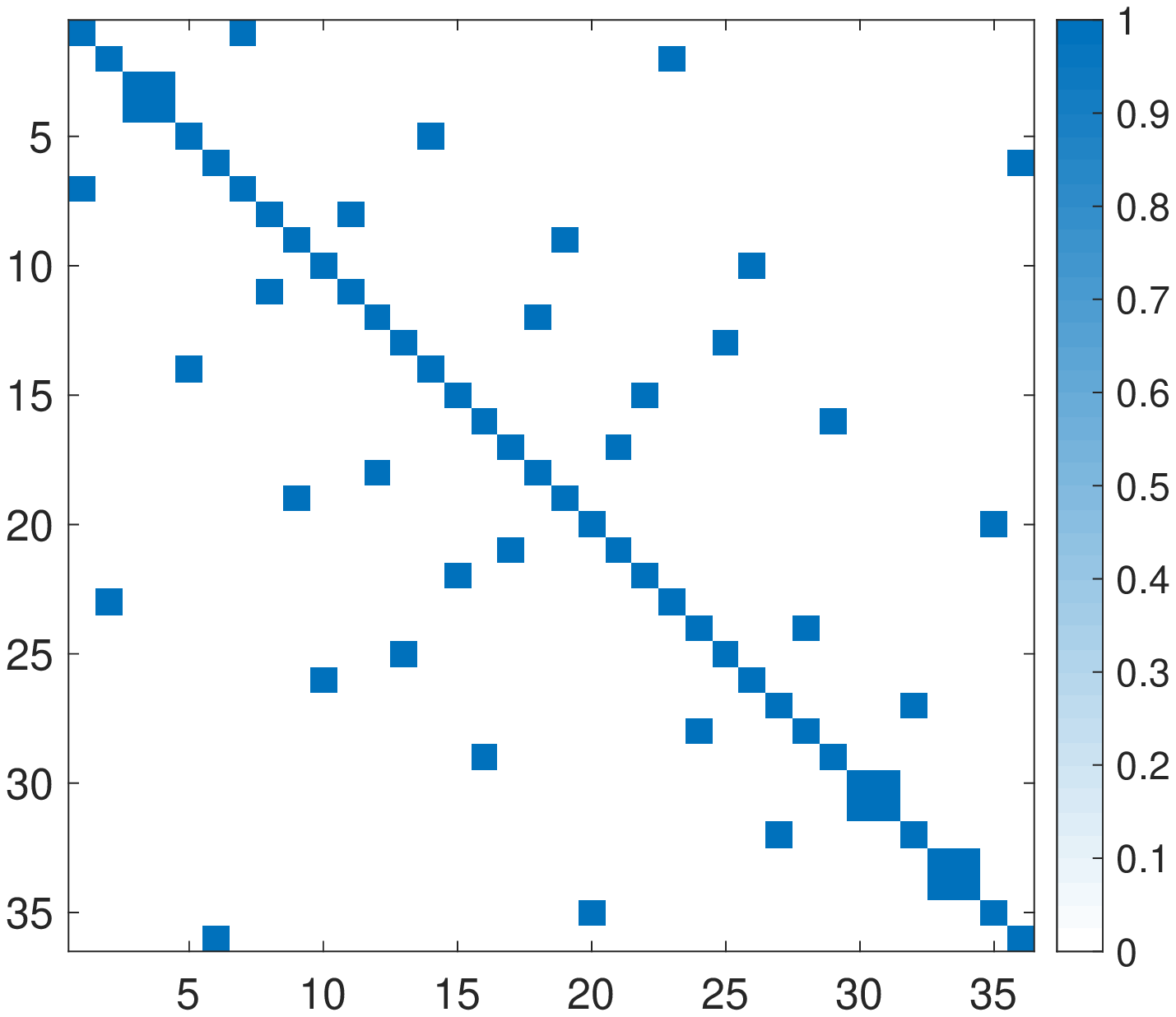}
    \includegraphics[width=0.32\textwidth]{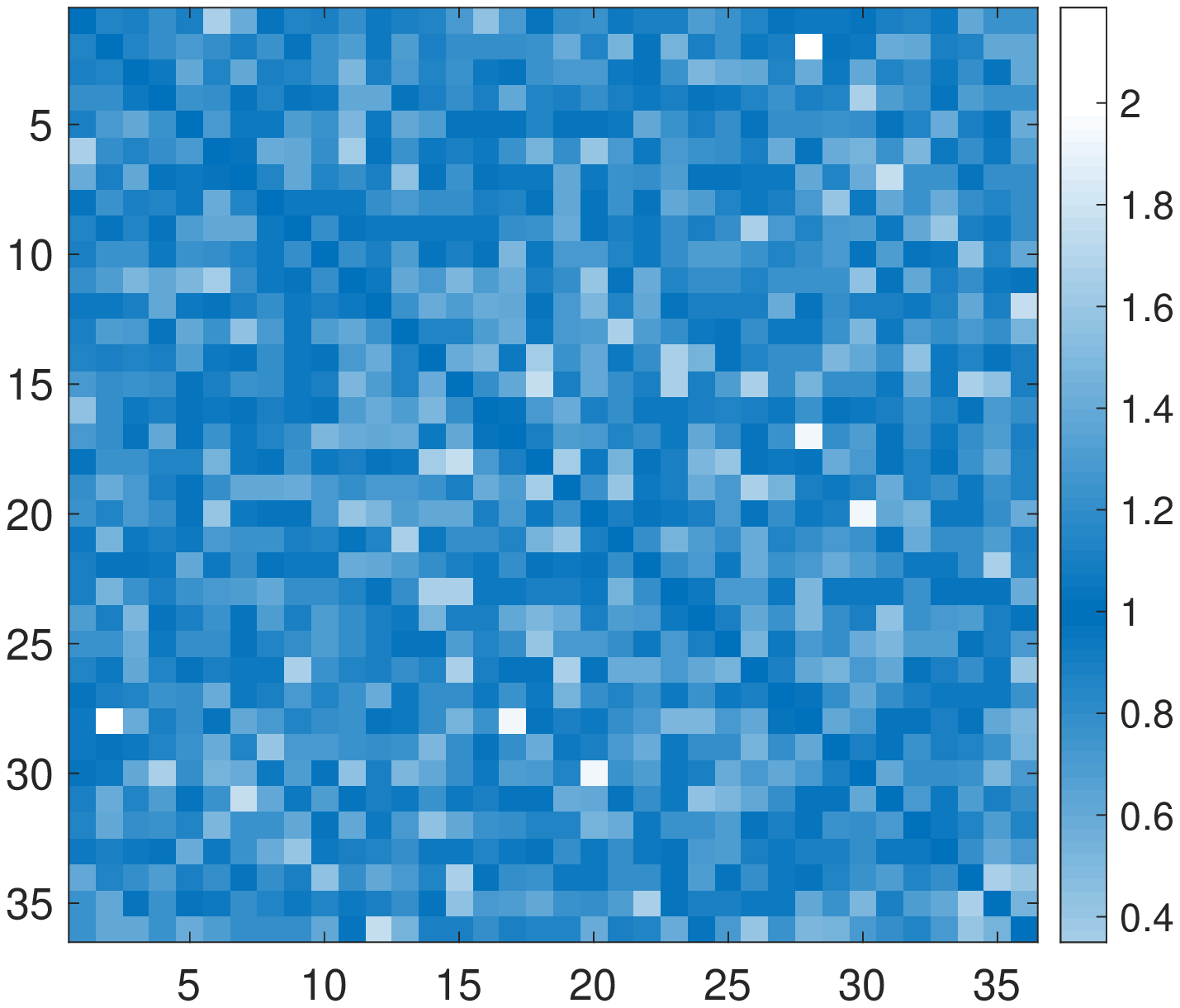}
  }
  \caption{Graphical representation of the adjacency matrix $A(t) = A_{ij}(t)$ for $36$ evaders with $P=2$ in the time intervals $[0,\Delta t]$ (left) and $[\Delta t,2\Delta t]$ (middle). The dotted points illustrate the interactions within the different batches, where the diagonal line is colored for reference. The right figure represents the averaged matrix over $400$ random batches.  The reduced model  in each time subinterval has a sparse network, while the complete system consists of all-to-all interactions.}
  \label{fig:velocity_field}
\end{figure*}

The reduced system has a network structure which switches the connectivity at each time $t_n$. In each time subinterval, $[t_n,t_{n+1}]$ the adjacency matrix is as follows,
\begin{equation}\label{eq:network}
\begin{aligned}
&A^n_{ij} := 
\begin{cases}
1 \quad \text{if }~ i \neq j ~\text{and}~ \text{in the same set among } {\mathcal C}_n^1, \ldots, {\mathcal C}_n^{N/P},\\
0 \quad \text{otherwise},
\end{cases}\\
&A_{ij}(t) := A^n_{ij} \quad \text{for }~ t \in [t_n,t_{n+1}),
\end{aligned}
\end{equation}
which is a piecewise constant matrix.

With this notation, the system \eqref{eq:sGBR} can be understood as a deterministic switching network system for $t \in [0,T]$,
\begin{equation}\label{eq:sGBR2} 
\begin{cases}
\dot{\mathbf x}_i = {\mathbf v}_i, \qquad i =1,\ldots,N,\\
\dot{\mathbf v}_i = \displaystyle\frac{1}{P-1}\sum_{k =1}^N A_{ki}(t) a({\mathbf x}_k-{\mathbf x}_i)({\mathbf v}_k-{\mathbf v}_i)
 + \frac{1}{P-1}\sum_{k=1}^N A_{ki}(t) g({\mathbf x}_k-{\mathbf x}_i)({\mathbf x}_k-{\mathbf x}_i)\\
\hspace{0.7em}\displaystyle - \frac{1}{M}\sum_{j=1}^M f({\mathbf y}_j-{\mathbf x}_i)({\mathbf y}_j-{\mathbf x}_i), \qquad i =1,\ldots,N,\\
\dot {\mathbf y}_j = {\mathbf u}_j(t), \qquad j=1,\ldots,M,\\
{\mathbf x}_i(t_n) = {\mathbf x}_i^n,\quad {\mathbf v}_i(t_n) = {\mathbf v}_i^n, \quad {\mathbf y}_j(t_n) = {\mathbf y}_j^n \qquad i =1,\ldots,N,\quad j=1,\ldots,M.
\end{cases}
\end{equation}
This system is well-posed, solutions exist and they are unique. Typically the velocity ${\mathbf v}$ is Lipschitz but not $C^1$ because of the discontinuities of the adjacency matrices.

The reduced system \eqref{eq:sGBR2}, including both the interactions among evaders and with drivers, involves a total number of  $O(N(P+M))$ interactions while, in the complete dynamics \eqref{eq:GBR}, the total number of interactions is  $O(N(N+M))$ per time step. 

Fig. \ref{fig:velocity_field} shows the adjacency matrix in \eqref{eq:sGBR2} among $36$ evaders, out of a choice of random batches of size $P=2$ (left and middle). The network structure changes at each time step, while the number of interactions is fixed to $N(P-1)$. When operating the numerical simulation with $T=4$ and $\Delta t = 0.01$, the connectivity along time $[0,4]$ is averaged $400$ times (right). As $\Delta t$ decreases to zero, it converges to the matrix of ones except for the diagonal elements, which is the adjacency matrix of the original system \eqref{eq:GBR}.

\subsection{The RBM for the optimal control problem}\label{sec:2-3}

As we mentioned in Section \ref{sec:2-1}, we are interested in finding the optimal control of the guiding problem. 
We adopt the RBM system above as the reduced model and we apply the gradient descent method. 
We then need to compute the gradient of the cost function \eqref{eq:cost} on the reduced RBM model \eqref{eq:sGBR}. The gradient can be derived, following the Pontryagin maximum principle \cite{pontryagin2018mathematical,trelat2005controle}, from the time evolution in the adjoint system of \eqref{eq:sGBR}.

First, we consider the adjoint of the original system \eqref{eq:GBR}. From the Pontryagin maximum principle, if we denote the controlled system \eqref{eq:GBR} and the cost \eqref{eq:cost} as
\begin{equation}\label{eq:PMP}
\begin{aligned}
\begin{cases}
\dot {\mathbf x} = F^{\mathbf x} ({\mathbf x},{\mathbf v},{\mathbf y}),\\
\dot {\mathbf v} = F^{\mathbf v} ({\mathbf x},{\mathbf v},{\mathbf y}),\\
\dot {\mathbf y} = F^{\mathbf y} ({\mathbf x},{\mathbf v},{\mathbf y}),\end{cases}
\text{and}\quad J = \int_0^T L({\mathbf x},{\mathbf v},{\mathbf y})dt,
\end{aligned}
\end{equation}
then its formal adjoint system can be described by the backward equations for $t \in [0,T]$,
\begin{equation}\label{eq:adjoint}
\begin{aligned}
\begin{cases}
- \dot {\mathbf p}^\mathsf{T} = {\mathbf p}^\mathsf{T} \nabla_{{\mathbf x}} F^{\mathbf x} + {\mathbf q}^\mathsf{T} \nabla_{{\mathbf x}} F^{\mathbf v} + {\mathbf r}^\mathsf{T} \nabla_{{\mathbf x}} F^{\mathbf y} + \nabla_{{\mathbf x}} L,\\
- \dot {\mathbf q}^\mathsf{T} = {\mathbf p}^\mathsf{T} \nabla_{{\mathbf v}} F^{\mathbf x} + {\mathbf q}^\mathsf{T} \nabla_{{\mathbf v}} F^{\mathbf v} + {\mathbf r}^\mathsf{T} \nabla_{{\mathbf v}} F^{\mathbf y} + \nabla_{{\mathbf v}} L,\\
- \dot {\mathbf r}^\mathsf{T} = {\mathbf p}^\mathsf{T} \nabla_{{\mathbf y}} F^{\mathbf x} + {\mathbf q}^\mathsf{T} \nabla_{{\mathbf y}} F^{\mathbf v} + {\mathbf r}^\mathsf{T} \nabla_{{\mathbf y}} F^{\mathbf y} + \nabla_{{\mathbf y}} L,\\
{\mathbf p}^\mathsf{T}(T) = 0,~{\mathbf q}^\mathsf{T}(T) = 0,~{\mathbf r}^\mathsf{T}(T) = 0,
\end{cases}
\end{aligned}
\end{equation}
where the variables ${\mathbf p}(t) \in \mathbb R^{Nd}$, ${\mathbf q}(t) \in \mathbb R^{Nd}$ and ${\mathbf r}(t) \in \mathbb R^{Md}$ are the adjoint states of ${\mathbf x}$, ${\mathbf v}$ and ${\mathbf y}$, respectively, and ${\mathbf p}^\mathsf{T}$ denotes the transpose of ${\mathbf p}$.

In our case,  the adjoint system of the complete dynamics \eqref{eq:GBR} would take the  following form in the dual variables $({\mathbf p}, {\mathbf q}, {\mathbf r})$ corresponding to $({\mathbf x}, {\mathbf v}, {\mathbf y})$, for each $i=1,\ldots,N$ and $j=1,\ldots,M$ :
\begin{equation*}
\begin{aligned}
-\dot{\mathbf p}_i^\mathsf{T} &= \frac{1}{N-1}\sum_{k =1,k\neq i}^N {\mathbf q}^\mathsf{T}_k 
\Big[ 
\nabla_{{\mathbf x}_i} a({\mathbf x}_k - {\mathbf x}_i)({\mathbf v}_k - {\mathbf v}_i))
+ \nabla_{{\mathbf x}_i} ( g({\mathbf x}_k - {\mathbf x}_i)({\mathbf x}_k - {\mathbf x}_i) ) \Big]
\\
&\quad - \frac{1}{M}  \sum_{j =1}^M
{\mathbf q}_i^\mathsf{T} \nabla_{{\mathbf x}_i} ( f({\mathbf y}_j - {\mathbf x}_i)({\mathbf y}_j - {\mathbf x}_i) )
+ \frac{2\alpha_1}{N}({\mathbf x}_i-{\mathbf x}_f), \qquad i =1,\ldots,N,\\
-\dot{\mathbf q}_i^\mathsf{T} &= {\mathbf p}_i^\mathsf{T} + \frac{1}{N-1}\sum_{k =1,k\neq i}^N {\mathbf q}^\mathsf{T}_k a({\mathbf x}_i - {\mathbf x}_k), \qquad i =1,\ldots,N,\\
-\dot{\mathbf r}_j^\mathsf{T} &= -\frac{1}{M}\sum_{k =1}^N {\mathbf q}^\mathsf{T}_k \nabla_{{\mathbf y}_j} ( f({\mathbf y}_j - {\mathbf x}_i)({\mathbf y}_j - {\mathbf x}_i) )
+ \frac{2\alpha_3}{M}({\mathbf y}_j-{\mathbf x}_f), \qquad j =1,\ldots,M.
\end{aligned}
\end{equation*}

Finally, the gradient of the cost function \eqref{eq:cost} to the control function ${\mathbf u}(t)$ is derived by
\begin{equation}\label{eq:gradient}
\begin{aligned}
\displaystyle\nabla_{\mathbf u} J &= \nabla_{\mathbf u}[{\mathbf p} \cdot F^{\mathbf x} +{\mathbf q} \cdot F^{\mathbf v} +{\mathbf r} \cdot F^{\mathbf y} + L]
= {\mathbf r} + (\alpha_2/M){\mathbf u}.
\end{aligned}
\end{equation}

Therefore, the implementation of the gradient descent method for the  complete system \eqref{eq:GBR} would lead to an iterative approximation of the optimal control as below
\begin{equation}\label{eq:GD}
\begin{aligned}
{\mathbf u}^{k+1} := {\mathbf u}^k- \alpha \nabla_{\mathbf u} J({\mathbf u}^k), \quad k \ge 0
\end{aligned}
\end{equation}
out of and initial guess ${\mathbf u}^0$ (that we could take to be ${\mathbf u}^0=0$) and with $\alpha >0$ small enough.

But, as we have described in the Introduction, this algorithm is computationally expensive since, in each step of the gradient descent iteration, it requires to solve the full state equation and adjoint system taking account of all interactions. We now present the adaptation of the gradient descent to the RBM reduced model.

We start from the RBM approximation \eqref{eq:sGBR2} in Section \ref{sec:2-2}. The corresponding adjoint system is computed piecewise in each of the time subintervals $[t_n, t_{n+1}]$, reducing the interaction terms to those corresponding to each batch. The adjoint dynamics of the reduced system for each $i=1,\ldots,N$ then reads as follows:
\begin{equation*}
\begin{aligned}
-\dot{\mathbf p}_i^\mathsf{T} &= \frac{1}{P-1}\sum_{k =1}^N A_{ki}(t) {\mathbf q}^\mathsf{T}_k 
\Big[ 
\nabla_{{\mathbf x}_i} a({\mathbf x}_k - {\mathbf x}_i)({\mathbf v}_k - {\mathbf v}_i))
+ \nabla_{{\mathbf x}_i} ( g({\mathbf x}_k - {\mathbf x}_i)({\mathbf x}_k - {\mathbf x}_i) ) \Big]
\\
&\quad - \frac{1}{M}  \sum_{j =1}^M
{\mathbf q}_i^\mathsf{T} \nabla_{{\mathbf x}_i} ( f({\mathbf y}_j - {\mathbf x}_i)({\mathbf y}_j - {\mathbf x}_i) )
+ \frac{2\alpha_1}{N}({\mathbf x}_i-{\mathbf x}_f), \qquad i =1,\ldots,N,\\
-\dot{\mathbf q}_i^\mathsf{T} &= {\mathbf p}_i^\mathsf{T} + \frac{1}{P-1}\sum_{k =1}^N A_{ki}(t) {\mathbf q}^\mathsf{T}_k a({\mathbf x}_i - {\mathbf x}_k), \qquad i =1,\ldots,N,
\end{aligned}
\end{equation*}
 the equation for ${\mathbf r}$ remaining the same.

\begin{remark}\label{r:RBM2}
Note that the adjoint system of the reduced RBM model coincides with the reduced RBM model of the complete adjoint system with the same adjacency matrix $A(t)$.
\end{remark}

 As in the forward dynamics, the RBM reduces the computational cost to $O(N(P+M))$ in the computation of the adjoint system. The gradient descent iteration can then be defined similarly as well in order to find the optimal control of the reduced model \eqref{eq:sGBR}.

However, since the approximation error accumulates in time, the resulting control may not suffice to guide the complete system \eqref{eq:GBR} in a long-time horizon.
The MPC procedure is now presented so to deal with the errors induced by the difference of the complete and reduced dynamics.

\subsection{The MPC procedure for the approximated model}\label{sec:2-4}

MPC is aimed to adapt the control  obtained for the reduced dynamics \eqref{eq:sGBR} to the full system \eqref{eq:GBR} in an iterative manner.

For this to be done, first, we choose a short time duration $\tau$ which determines the length of the control time intervals in which the reduced control will be applied to the complete system. But these controls are computed on longer time intervals of length $\hat T \geq \tau$. 
The time parameters $\tau$ and $\hat T$ critically affect the performance of the MPC strategy though there is no general argument to determine them. On one hand, to get reliable control functions, we need to set a short time $\tau$ and a long time $\hat T$. For example, in the numerical simulations of Section \ref{sec:simul}, we set $T = 4$, $\tau = 1.5$ and $\hat T = 3$. 

With these choices of the parameters $\tau$ and $\hat T$, we define the discrete times $\tau_m := m\tau$ for $m \geq 0$. In the first time interval $[\tau_0=0, \tau_1=\tau]$, we apply the control computed in $[0, \hat T]$ for the RBM reduced dynamics \eqref{eq:sGBR} into the complete one \eqref{eq:GBR}. This leads to a given state of the complete system \eqref{eq:GBR} in time $\tau_1$. Taking that one as initial state, we recompute the control of the RBM system in the time horizon $[\tau_1, \tau_1+ \hat T]$. And this strategy is repeated again and again for $m\ge 0$ until $\tau_m \geq T$. 

In this manner, the MPC strategy takes account of the gap between the trajectories of the complete and the reduced dynamics.

Note that a proper choice of $\tau$ depends on the errors on the reduced dynamics and the control functions. If the reduced model gets more accurate, then we may choose a bigger $\tau$ so that the computation becomes cheaper.
On the other hand, $\hat T$ is rather affected by the nature of the system \eqref{eq:GBR}. If $\hat T$ is too short, the control would have a big difference from the optimal one, unable to capture the global requirements of the control problem. 

\subsection{Implementation of the MPC-RBM algorithm}
\label{sec:2-5}

We now summarize the discussion from Section \ref{sec:2-1} to Section \ref{sec:2-4} on the RBM-MPC algorithm: 

\begin{enumerate}
\item
The objective of the algorithm is to control the system \eqref{eq:GBR} with a small computational cost, solving the optimal control problem for  the cost $J({\mathbf u})$ of \eqref{eq:cost}  in a (possibly long) time-interval $[0,T]$. We do it combining MPC with the RBM approximated dynamics.
\item
We set a short time length $\tau>0$ for the control time and a long time length $\hat T>0$ for the predictive time. Then, we define the discrete times $\tau_m := m\tau$, $m \ge 0$. For each $m$, we iteratively solve the optimal control problems on each predictive time interval  as follows.
\item
The original dynamics  being autonomous in each step of the iteration the time-interval $[\tau_m,\tau_m+\hat T]$ can be shifted to  $[0,\hat T]$. We adopt the control of the RBM approximation as initial guess of the control ${\mathbf u}(t)$.

The RBM model is implemented as described above, combined  with the Euler forward time-discretization.

The RBM control is computed by a gradient descent method that needs of  corresponding RBM adjoint system that we described above.

\item This control is implemented in the com plate system up to time $\tau_1$. 

\item
This strategy is iteratively applied in the intervals $[\tau_m,\tau_{m}+\hat T]$ for the whole time interval $[0,T]$.

\end{enumerate}

This Algorithm \ref{alg:MPC-RBM} has the following features:

\begin{enumerate}
\item 
The cost function of the recursive minimization of $J$ monotonically decreases during the optimization process since we use the standard gradient descent algorithm as in \cite{nesterov2013introductory} to the RBM reduced dynamics.
\item
The resulting control is not the (local) minimizer of the cost $J$ for the complete dynamics. But it gives an approximative.
\item
The total computational effort is reduced by a factor $O(N(P+M))/O(N(N+M))$, $N$ being the dimension of the complete system and $P$ the size of the batches.
\item 
The process of MPC acts in a semi-feedback way, where we periodically check the current states of the system to update the control. 
This strategy is  adaptive so that it could cope with unexpected uncertain events in the global dynamics.

\end{enumerate}

\begin{algorithm}
\begin{algorithmic}
 \caption{MPC-RBM algorithm for \eqref{eq:GBR}}
\State $P$, $\tau$, $\Delta t$ and $\hat T$ are given.
\Function {RBM-State}{${\mathbf x}^0,{\mathbf v}^0,{\mathbf y}^0,{\mathbf u}(t)$}
  \State Fix a random seed for the choices of batches. 
  \For {$n$ from $0$ to $[\hat T/\Delta t]$} %
    \State Divide $\{1,2,\ldots,N\}$ into random batches with size $P$.
    \For{each batch} %
      \State Update ${\mathbf x}_i$ and ${\mathbf v}_i$ by solving the reduced model \eqref{eq:sGBR} from $t = n \Delta t$ to $t= (n+1) \Delta t$.
    \EndFor
    \State Update ${\mathbf y}_i$ from $t = n \Delta t$ to $t= (n+1) \Delta t$.
  \EndFor
  \State \Return ${\mathbf x}(t)$, ${\mathbf v}(t)$ and ${\mathbf y}(t)$.
\EndFunction
\State

\Function {RBM-Costate}{${\mathbf p}^0,{\mathbf q}^0,{\mathbf r}^0,{\mathbf x}(t),{\mathbf v}(t),{\mathbf y}(t)$)}
  \State Fix a random seed, the same one from RBM-State in reverse order.
  \For{$n$ from $[\hat T/\hat \tau]$ to $0$} %
    \State Divide $\{1,2,\ldots,N\}$ into random batches with size $P$.
    \For{each batch} %
      \State Update ${\mathbf p}_i$ and ${\mathbf q}_i$ (the adjoint of ${\mathbf x}_i$ and ${\mathbf v}_i$) by solving the adjoint of the reduced model \eqref{eq:sGBR} from $t= (n+1) \Delta t$ to $t = n \Delta t$.
    \EndFor
    \State Update ${\mathbf r}_i$ (the adjoint of ${\mathbf y}_i$) from $t= (n+1) \Delta t$ to $t = n \Delta t$.
  \EndFor
  \State \Return ${\mathbf p}(t)$, ${\mathbf q}(t)$ and ${\mathbf r}(t)$.
\EndFunction
\State

\Function {OCP}{${\mathbf x}^0,{\mathbf v}^0,{\mathbf y}^0,{\mathbf u}_0(t)$)}
  \State Define the cost function $J$ with \eqref{eq:cost} over $[0,\hat T]$.
  \State Initialize the control ${\mathbf u}(t)$ with a guess ${\mathbf u}_0(t)$.
    \While{$\|D_{{\mathbf u}} J\| < \varepsilon$ (or any stopping criteria)} %
      \State Operate RBM-State and RBM-Costate.
      \State Calculate the gradient $D_{{\mathbf u}} J$ of the cost $J$.
	    \State Update ${\mathbf u}(t)$ using $D_{{\mathbf u}} J$.
    \EndWhile
	\State \Return ${\mathbf u}(t)$.
\EndFunction
\State

\Procedure {MPC-RBM algorithm}{}
	\State Set the initial data ${\mathbf x}^0,{\mathbf v}^0,{\mathbf y}^0$ for the system \eqref{eq:GBR}.
	\State Give an initial guess on control ${\mathbf u}_0(t)$ for $t \in [0,T]$.
	\State Let $\tau_m := m \tau$ for $m=0,1,\ldots,[T/\tau]+1$.
  \For{$m$ from $0$ to $[T/\tau]$} %
	  \State Operate OCP with {${\mathbf x}(\tau_m),{\mathbf v}(\tau_m),{\mathbf y}(\tau_m),{\mathbf u}_0(t)$} to get the optimal control ${\mathbf u}(t)$ for $t \in [\tau_m,\tau_m+\hat T]$.
    \State Process the original system \eqref{eq:GBR} with the control ${\mathbf u}(t)$ for $t \in [\tau_m,\tau_{m+1}]$.
  \EndFor
  \State \Return the trajectories and control over $t \in [0,T]$.
\EndProcedure
\label{alg:MPC-RBM}
\end{algorithmic}
\end{algorithm}
 
\subsection{An open problem: error analysis of the MPC-RBM method}\label{sec:2-6}

Developing a detailed error analysis of the MPC-RBM as a function of the various parameters entering in the process if is a challenging issue beyond the scope of this article.

Indeed, in \cite{jin_random_2020}, the error analysis of the RBM is developed under suitable contractility properties of the original dynamics. This analysis does not apply directly to the model under consideration.

In the numerical simulations \cite{carrillo_consensus-based_2019,ha_convergence_2019,jin_random_2020}, the RBM shows a successful guess on the long-time behavior of the original systems. Hence, the uniform-in-time error analysis of the RBM is an interesting open problem on collective behavior models, for example, the aggregation dynamics or pedestrian flows.

In addition to that, one should further consider the error induced by the MPC strategy. It needs a sensitivity analysis from the trajectory error to the control function, which has been studied for linear problems \cite{mhaskar2006robust,prett1987design}. 

In the absence of analytical results on the MPC-RBM algorithm, we present several computational experiments that confirm the numerical efficiency of the method in the following section.

\section{Simulations on the MPC-RBM algorithm}\label{sec:simul}

We consider the guiding problem with $36$ evaders ($N=36$) and $2$ drivers ($M=2$) in the two-dimensional space ($d=2$). The target is chosen to be  ${\mathbf x}_f = (0.5,0.5)$ with the final time $T=4$ and the time step $\Delta t = 0.01$ in the implementation of the RBM. Initially the evaders are uniformly distributed in $[-0.2,0.2]^2$ with zero velocities, and the drivers start from two points $(-1,0)$ and $(0,-1)$. 

The cost function is given by \eqref{eq:cost} with the regularization coefficients $\alpha_1 = 1$ and $\alpha_2=\alpha_3 = 10^{-4}$ as follows:
\begin{equation*}
\begin{aligned}
J({\mathbf u}) :=& \int_0^T\left[ \frac{1}{N}\sum_{k=1}^N |{\mathbf x}_k-{\mathbf x}_f|^2
+ \frac{10^{-4}}{M}\sum_{j=1}^M |{\mathbf u}_j|^2 + \frac{10^{-4}}{M}\sum_{j=1}^M |{\mathbf y}_j-{\mathbf x}_f|^2 \right]dt. 
\end{aligned}
\end{equation*}

The numerical simulations are operated in Matlab with a laptop consisting of CPU i5-4258U 2.4GHz and RAM DDR3L 8GB 1600MHz. The symbolic calculations on gradients and adjoints are implemented with CasADi \cite{Andersson2018}. The functions on the forward and adjoint dynamics (\emph{RBM-STATE} and \emph{RBM-COSTATE} in Algorithm \ref{alg:MPC-RBM}) are also implemented as CasADi symbolic functions for a fast calculation in a pre-calculated form. The random batches are chosen with the \emph{randperm} function in Matlab.

\subsection{Simulations on RBM for the controlled dynamics}\label{sec:3-1}

First, we explore the performance of the RBM for various values of $P$ in the controlled dynamics, to discuss its time efficiency and approximation error.

\begin{figure*}
  \centering
  {
    \includegraphics[width=0.45\textwidth]{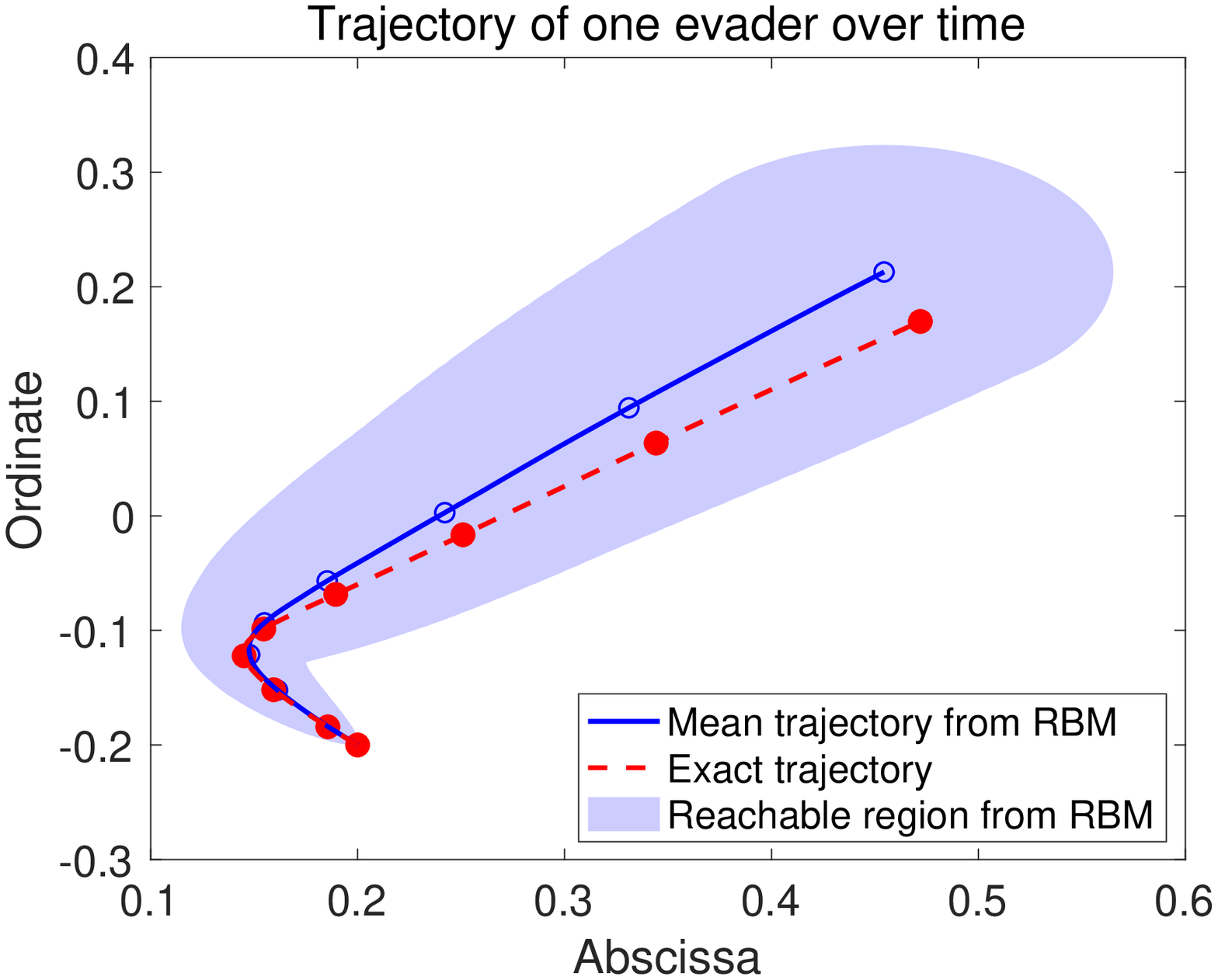}
    \includegraphics[width=0.45\textwidth]{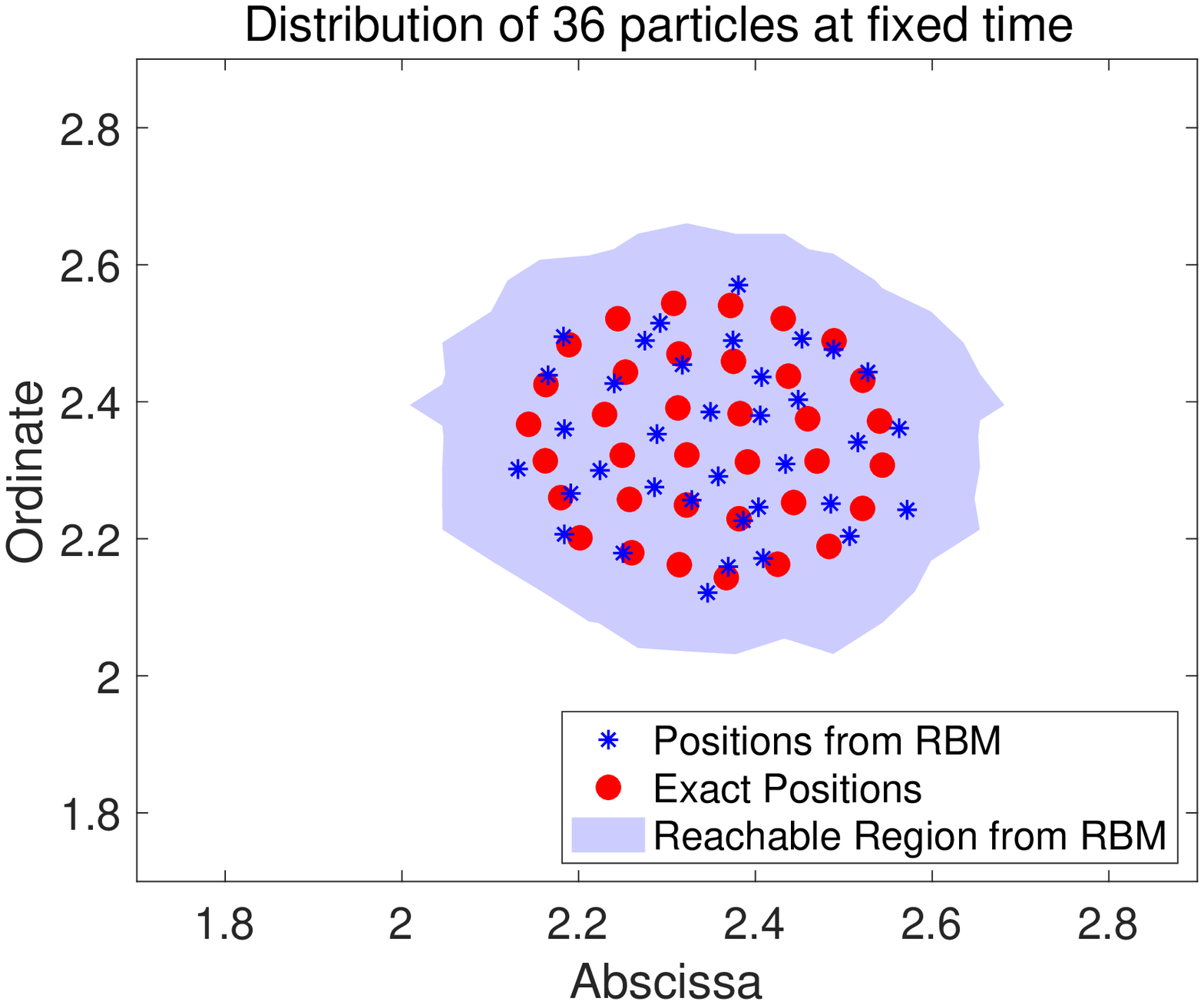}
  }
  \caption{A simulation on trajectories of the evaders using the RBM ($P=2$). \emph{Left}: the blue line shows the mean trajectory of one evader starting from $(-0.2,0.2)$ for $t\in[0,4]$, which is averaged among $200$ RBM approximations. The red line represents the trajectory from \eqref{eq:GBR}. The reachable region is drawn with $95\%$ reliability from the standard deviation and the 2D normal distribution. \emph{Right}: The blue marks show the final positions of the $36$ evaders from the reduced RBM model at $t=10$, while the red marks are the final positions from the original system. The reachable region draws the area containing more than one particle from $200$ RBM approximations, which represents $99\%$ credible region.\\}
  \label{fig:error1}
\end{figure*}

\begin{figure*}
  \centering
  {
    \includegraphics[width=0.45\textwidth]{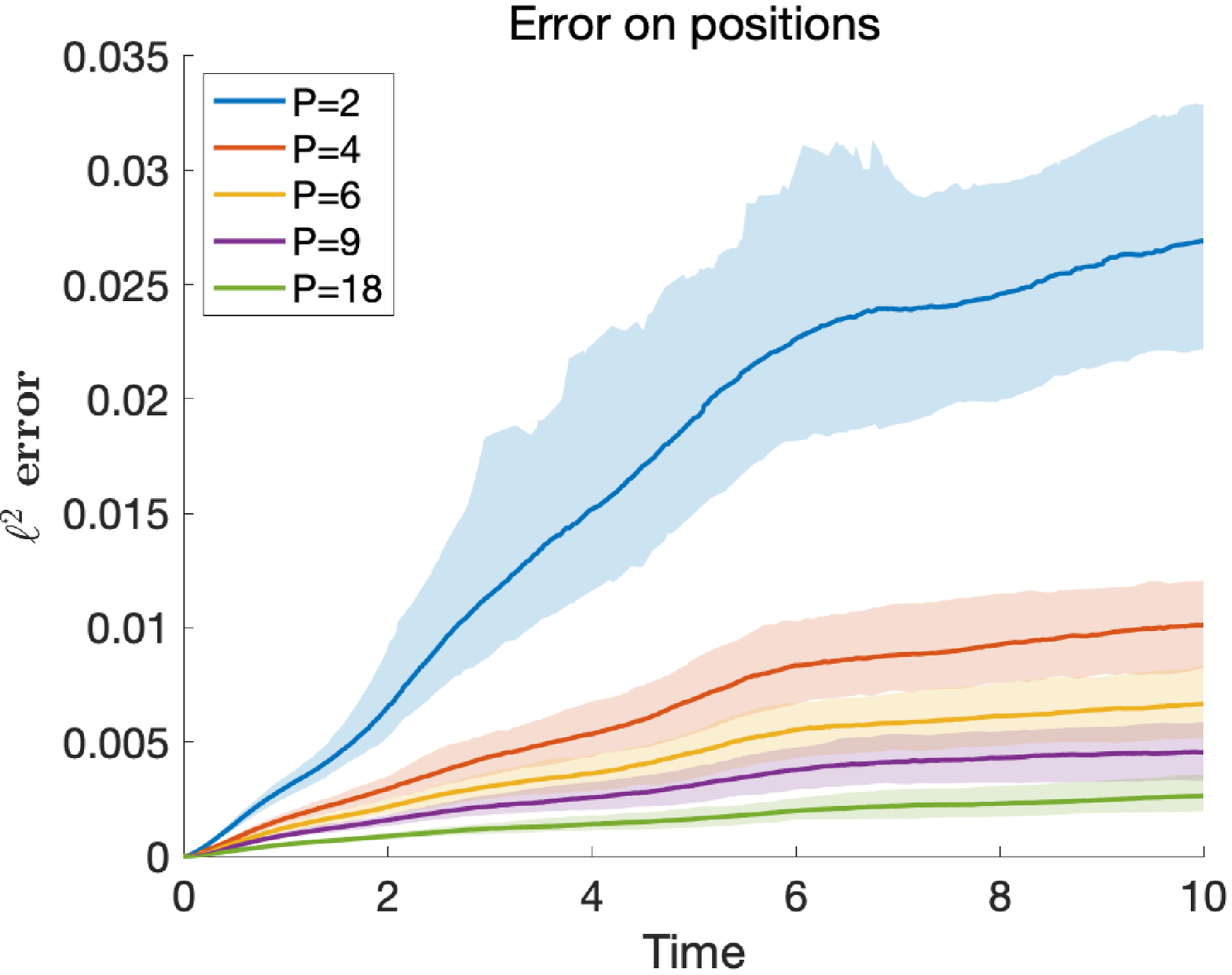}
    \includegraphics[width=0.45\textwidth]{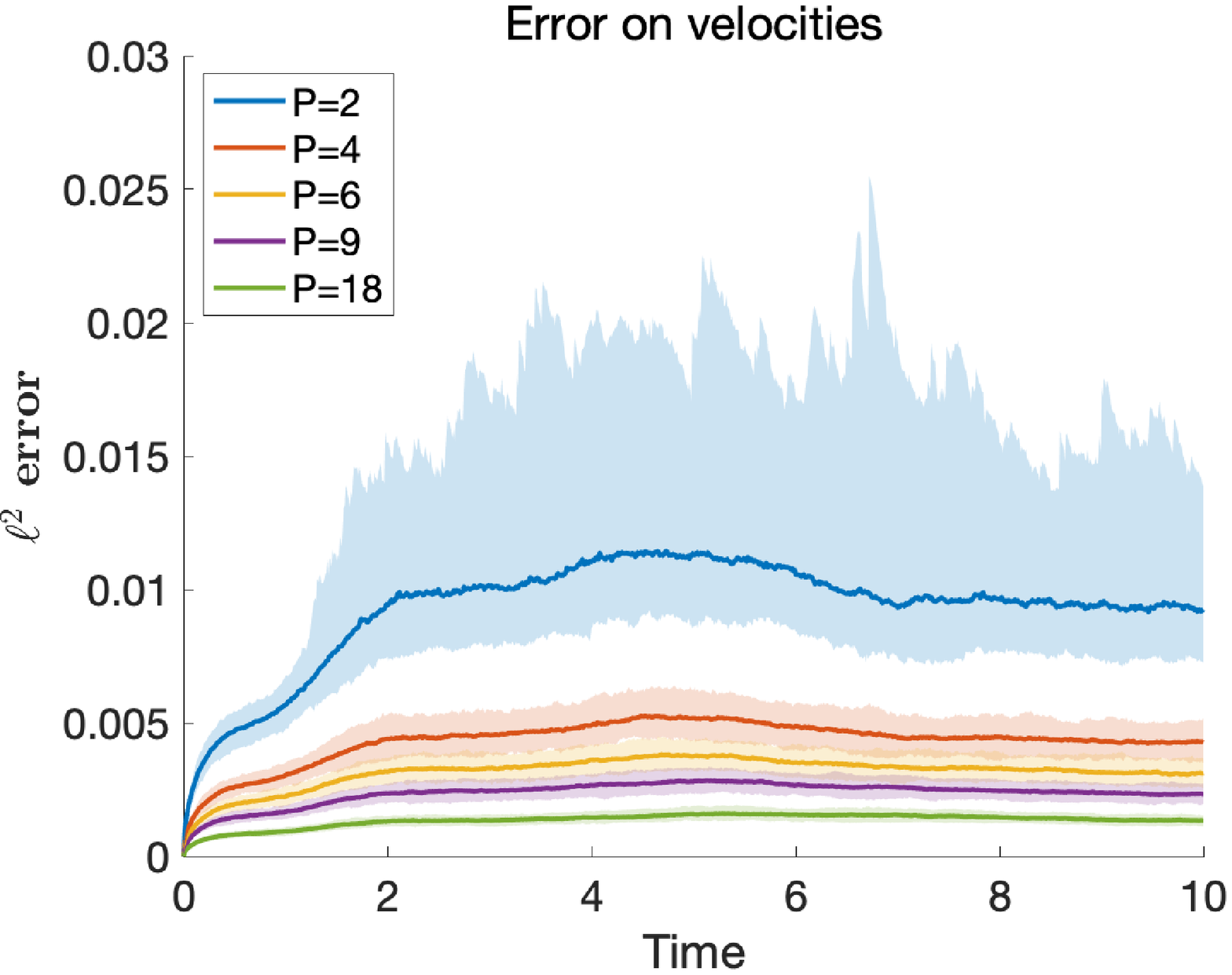}
  }
  \caption{The comparison between the dynamics from the RBM approximations and the exact system over time $t \in [0,10]$. \emph{Left}: the $\ell^2$ errors on positions at each time for various $P$ in RBM, where the median values over $200$ random simulations drawn with lines. The colored region shows the confidence interval for each time with $95\%$ reliability. \emph{Right}: the corresponding $\ell^2$ errors for velocities.\\}
  \label{fig:error2}
\end{figure*}

\begin{table*}
\begin{center}
\begin{tabular}{l|*{4}{c}}
Time-evolution    & Computation time & (Ratio) & Interactions & (Ratio) \\
\hline
Full system & 37.7220 & (9.759) & 8136 & (10.27) \\
RBM ($P=2$) & \phantom{0}3.8654 & (1.000) & \phantom{0}792 & (1.000) \\
RBM ($P=4$) & \phantom{0}6.6993 & (1.733) & 1224 & (1.545) \\
RBM ($P=6$) & \phantom{0}8.9043 & (2.304) & 1656 & (2.091) \\
RBM ($P=9$) & 11.3447 & (2.935) & 2304 & (2.909) \\
RBM ($P=18$) & 19.4503 & (5.032) & 4248 & (5.364) \\
\end{tabular}
\end{center}
\caption{The computation time (in milliseconds) to calculate the forward dynamics (controlled trajectories) for $t \in [0,10]$. For each reduced RBM model, it also describes the number of calculated interactions per time step to compare with the complete dynamics. The standard Euler forward method is used and averaged for $1000$ simulations.}\label{tab:time1}
\end{table*} 

Fig. \ref{fig:error1} shows the positions of evaders simulated along time (left) and at fixed time (right). The trajectories are calculated with the original system and the reduced RBM model with the same control functions. The controls are given by ${\mathbf u}_1(t)=(0.2,0.02)$ and ${\mathbf u}_2(t) = (0.02,0.2)$ for $t \in [0,10]$, where the drivers push the evaders toward the northeast direction. 

The left figure shows the trajectory of one evader for $t \in [0,4]$ starting from the point $(-0.2,0.2)$. The blue colored reachable region is the confidence region with radius $1.73\sigma$ ($95\%$ reliability from the 2D normal distribution) from the standard deviation $\sigma$ at each time. 
In the simulation, the radius of the reachable region at $t=4$ is nearly $0.1$, which is too rough compared to the distance between nearby evaders.
However, according to the right figure, the blue markers (evaders from the reduced model) have a similar distribution to the red dots (evaders from the original system) at fixed time, $t=10$. In particular, the reduced model produces good approximations for the mean position and the diameter of the herd, which are critical in the guiding problem.

On the other hand, Fig. \ref{fig:error2} shows the approximation errors while the system evolves in time $t \in [0,10]$. The averaged $\ell^2$ errors per each evader are presented from $200$ independent RBM approximations for each $P$. 
The median value of the errors is drawn in lines, and the colored region represents the confidence interval at each time. Note that the errors on velocities are bounded and the errors on positions grow linearly, not exponentially. 
We can also observe that the errors are reduced and more stable when $P$ gets bigger.

In Table \ref{tab:time1}, the computation time for the evolution is described with different $P$, averaged over $1000$ simulations.
For $N$ evaders and $M$ drivers in a $d$-dimensional space, the number of interactions in the derivatives $(\dot{\mathbf x},\dot{\mathbf v},\dot{\mathbf y})$ can be estimated by $Nd(1+Md+P(d+1))$.
Therefore, from the full system to the RBM one with $P=2$, the computational  cost is divided, roughly,  by a factor of $10$  ($8136$ in the original system and $792$ in the reduced model). Note that the computational time is nearly proportional to the number of interactions in Table \ref{tab:time1}.

\subsection{Simulations on RBM for the optimal controls}

We now compute the optimal controls both for the original system and the reduced RBM model. As described above, we fix the choice of random batches during the optimization process.

\begin{figure*}
  \centering
  {
    \includegraphics[width=0.8\textwidth]{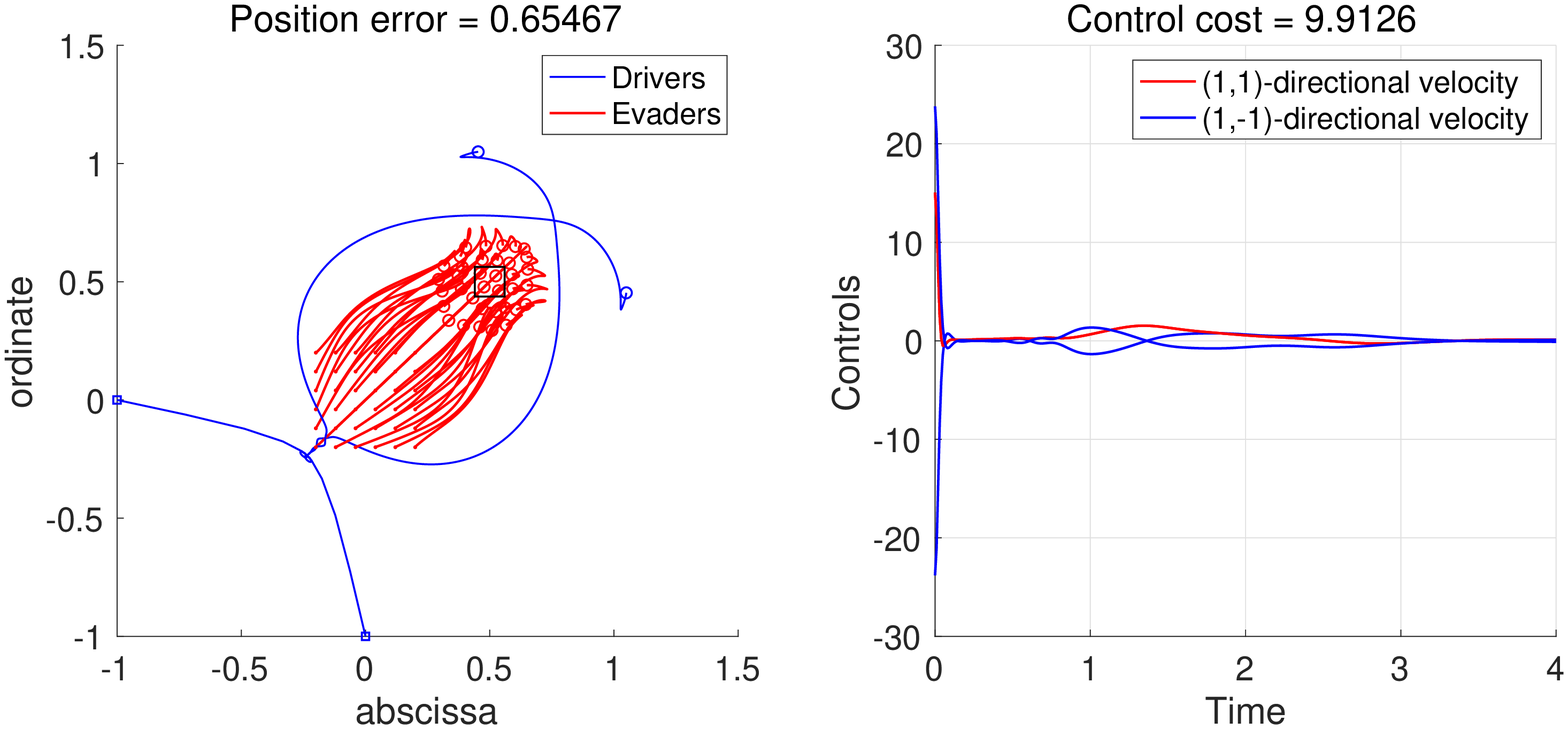}
    \includegraphics[width=0.8\textwidth]{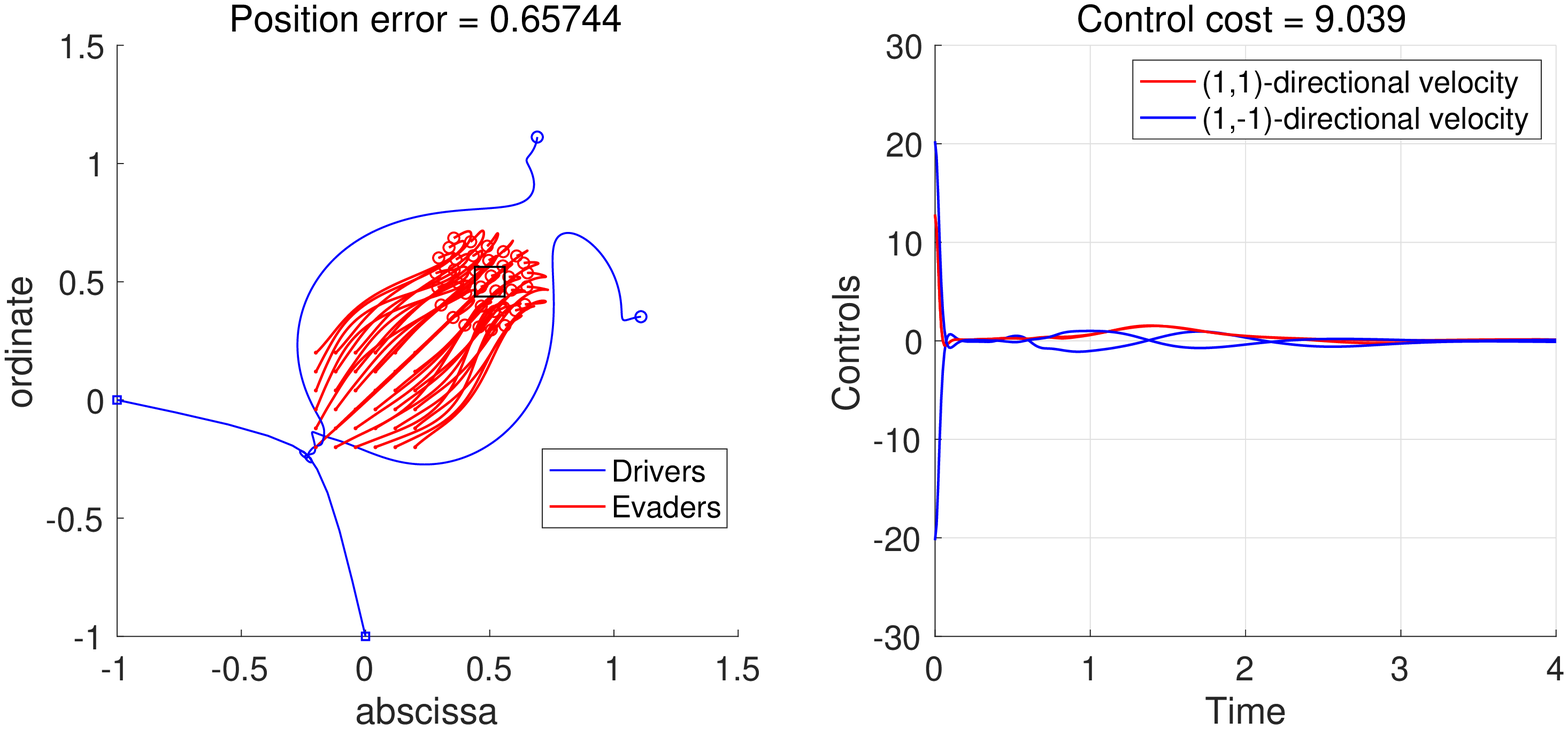}
  }
  \caption{The controlled trajectories from the optimal controls calculated with the original system (above) and the RBM model with $P=2$ (below). The trajectories are simulated on the original system. \emph{Left}: the trajectories of the drivers and evaders for $t \in [0,4]$. Red lines are for the evaders and blue lines are for the drivers. Initially the evaders are pushed in the northeast direction. The position error denotes the time integration of the averaged squared distance from the evaders to the target point. \emph{Right}: the control functions along time. Red lines represent the control in the direction of $(1,1)$ and blue lines are for its orthogonal direction. The control cost is the sum of time integrations, one on the averaged squared norm of the controls and the other on the averaged squared distance from the drivers to the target.}
  \label{fig:test1}
\end{figure*}

\begin{table*}
\begin{center}
\begin{tabular}{l|*{5}{c}}
Adaptive GD   & GD iterations & EV calculations & Cost $J$ & Computation time & (Ratio) \\
\hline
Full system  & 3123 & 9415 & 0.6646 & 767.96 & (8.004) \\
RBM ($P=2$) & 2721 & 7567 & 0.6665 & \phantom{0}95.95 & (1.000) \\
RBM ($P=4$) & 3054 & 9208 & 0.6651 & 164.02 & (1.710) \\
RBM ($P=6$) & 3027 & 9127 & 0.6651 & 214.05 & (2.231) \\
RBM ($P=9$) & 2604 & 7858 & 0.6648 & 229.57 & (2.393) \\
RBM ($P=18$) & 2654 & 8008 & 0.6648 & 392.71 & (4.093) \\
\end{tabular}
\end{center}
\caption{The computation time in the (adaptive) gradient descent algorithm for the guiding problem. The stopping criterion is set with tolerance $10^{-6}$ on the cost $J$, which tests the ratio between the difference of the cost and the previous cost. GD iterations are the number of iterations in gradient descent, and EV calculations are the number of calculations on the time evolution for the adaptive step size.}\label{tab:time2}
\end{table*}

In order to find the optimal control, here we use the gradient descent with an adaptive step size. The initial guess on the control is the constant functions used in the simulations of Fig. \ref{fig:error1}. The step size $\alpha$ in \eqref{eq:GD} is initially set to be $0.1$, and try a half of it if the next cost $J({\mathbf u}^{k+1})$ is not smaller than $J({\mathbf u}^{k})$ by the ratio of $10^{-6}$. For the next optimization step, we try the twice of the previous step size. This adaptive method guarantees the monotonicity of the cost $J$ during the optimization process. The algorithm will be terminated if $\alpha$ is less than $10^{-15}$.

Table \ref{tab:time2} shows the computation time to find the optimal controls with different values of $P$. Note that the optimization steps are not much different in all cases, hence, the computation time follows similar ratios to Table \ref{tab:time1}.

In Fig. \ref{fig:test1}, the controlled trajectories are presented with the different optimal controls from the original system and the reduced model ($P=2$). Both trajectories are simulated on the system \eqref{eq:GBR}. Note that the detailed motion of drivers differs, however, the overall running costs are similar, as presented in Table \ref{tab:time2}.

\subsection{The effect of MPC in the guiding problem}\label{sec:3-3}

\begin{figure*}
  \centering
  {
    \includegraphics[width=0.8\textwidth]{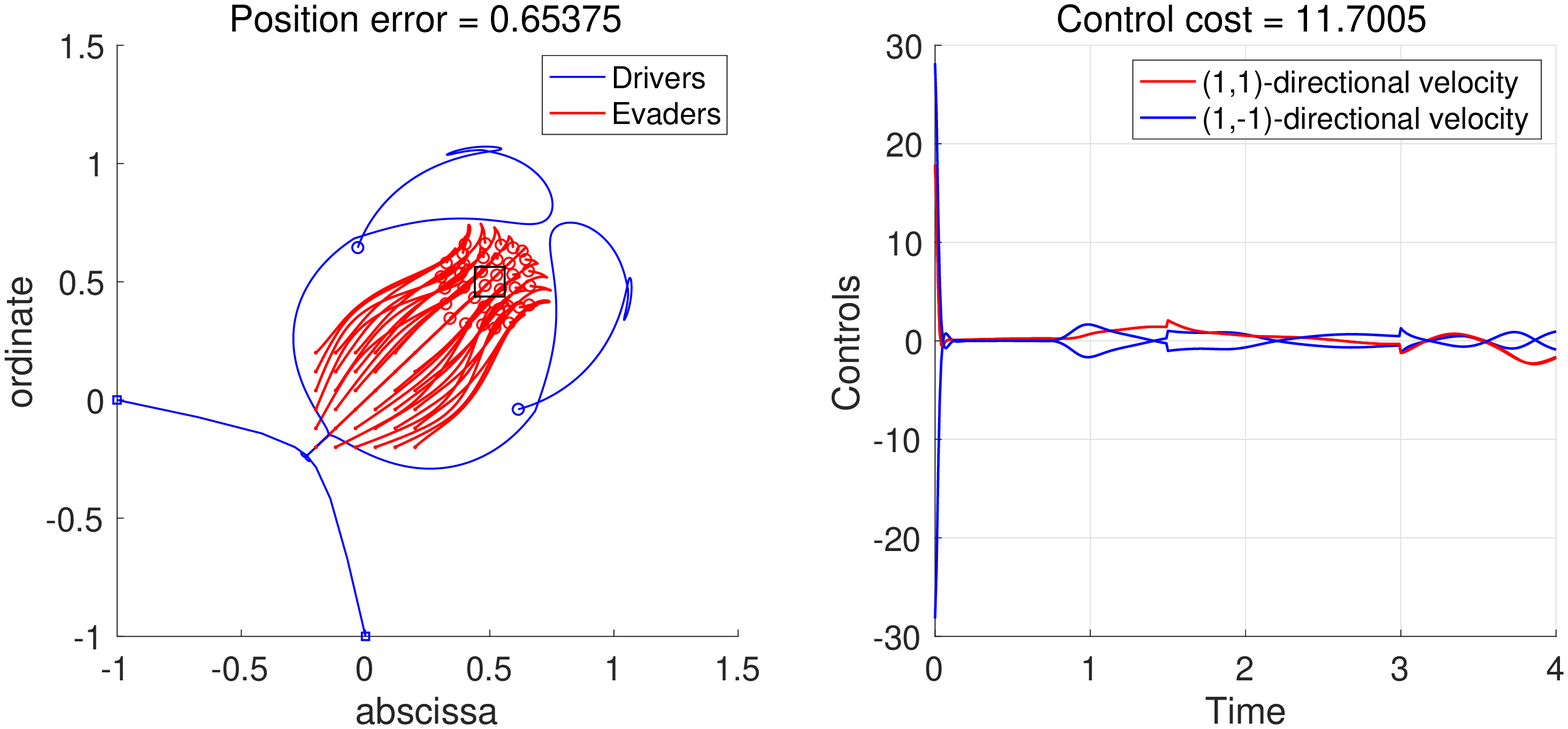}
    \includegraphics[width=0.8\textwidth]{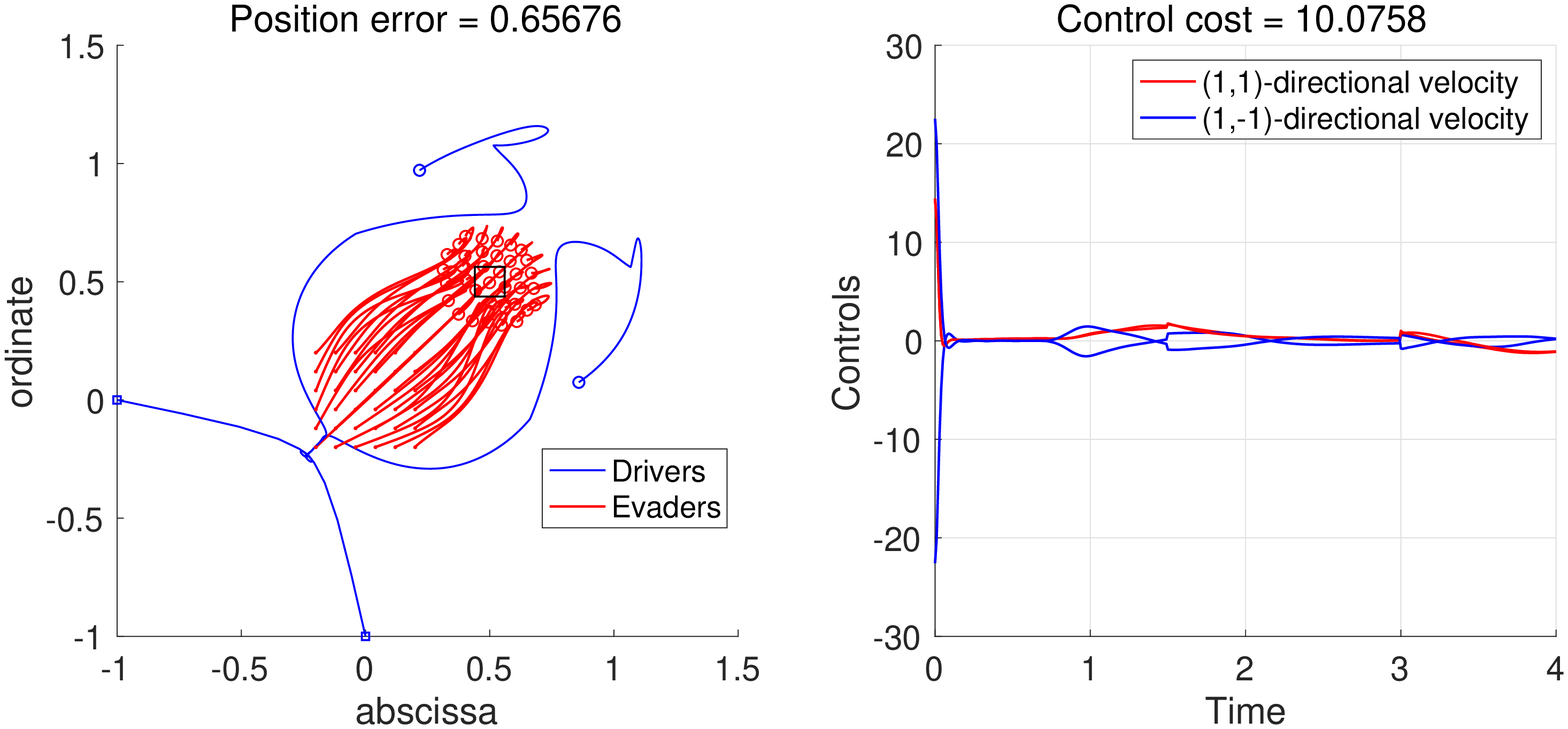}
  }
  \caption{The controlled trajectories from the MPC algorithms with $\tau = 1.5$ and $\hat T = 3$ over the original system (above) and the reduced RBM model with $P=2$ (below). \emph{Left}: the trajectories of the drivers and evaders for $t \in [0,4]$. \emph{Right}: the control functions along time.}
  \label{fig:test2}
\end{figure*}

The simulations in Fig. \ref{fig:test1} shows the open-loop optimal control, which does not use the strategy of MPC. Next, we compare the MPC-RBM algorithm with the same conditions. Fig. \ref{fig:test2} shows the simulations using the strategy of MPC, where the controls are calculated from the original system (above) and the RBM approximation with $P=2$ (below).

For the process of MPC, we set $\tau = 1.5$ and $\hat T = 3$. This implies that the controls on the intervals $[0.0,1.5]$, $[1.5,3.0]$ and $[3.0,4.0]$ are calculated from the predictive time intervals $[0.0,3.0]$, $[1.5,4.5]$ and $[3.0,6.0]$, respectively.
For the first interval $[0,3]$, we used the same initial guess on the control as before. However, for the next time intervals, we used the control obtained in the previous time intervals. In detail, for the second interval $[1.5,4.5]$, we used the optimal control from the time horizon $[0.0,3.0]$ to apply it on $[1.5,3.0]$, and set zero values on $[3.0,4.5]$.

One of the differences between Fig. \ref{fig:test1} and Fig. \ref{fig:test2} is the motion of the drivers near the final time. In Fig. \ref{fig:test2}, the drivers move around the evaders to capture them near the target point $(0.5,0.5)$. This is from the effect of $\hat T$ since the control is obtained with the dynamics on $t \in [0,6]$, not $[0,4]$. 

On the other hand, the computation time for Fig. \ref{fig:test2} is described in Table \ref{tab:time_mpc}. The MPC algorithm formulates the optimal control problems three times, but the overall times are similar to Table \ref{tab:time2}. In particular, the simulation of the MPC algorithm takes $956$ seconds with the original system and $83$ seconds with the RBM, while the open-loop controls are calculated with $768$ and $96$ seconds, respectively. 

Since it differs from case to case, it is difficult to estimate the computation time with MPC. However, note that the computational cost of the MPC-RBM algorithm is still in the order of $O(N(P+M))$, which is much less than $O(N^2)$ of the complete system when $N$ is large.

\begin{table*}[t]
\begin{center}
\begin{tabular}{l|*{3}{c}|c}
\hline
& Control interval & GD iterations & Computation time & Cost $J$ \\
\hline
MPC-Full system & $[0.0,1.5]$ & 2819 & 586.94 & 0.6654  \\
(Fig. \ref{fig:test2}, above) & $[1.5,3.0]$ & 1307 & 241.98 & \\
(Total time: $955.81$)               & $[3.0,4.0]$ & 7006 & 126.89 &  \\
\hline
MPC-RBM ($P=2$)       & $[0.0,1.5]$ & 1727 & \phantom{0}47.42 & 0.6668  \\
(Fig. \ref{fig:test2}, below) & $[1.5,3.0]$ & \phantom{0}216 & \phantom{00}7.30 & \\
(Total time: $82.96$)   & $[3.0,4.0]$ & 1040 & \phantom{0}28.24 & \\
\hline
\end{tabular}
\end{center}
\caption{The number of iterations and computation time (in seconds) for Fig. \ref{fig:test2} to find the optimal controls with MPC.}\label{tab:time_mpc}
\end{table*}

\subsection{Simulations of MPC-RBM on a noisy system}\label{sec:3-4}

\begin{figure*}
  \centering
  {
    \includegraphics[width=0.8\textwidth]{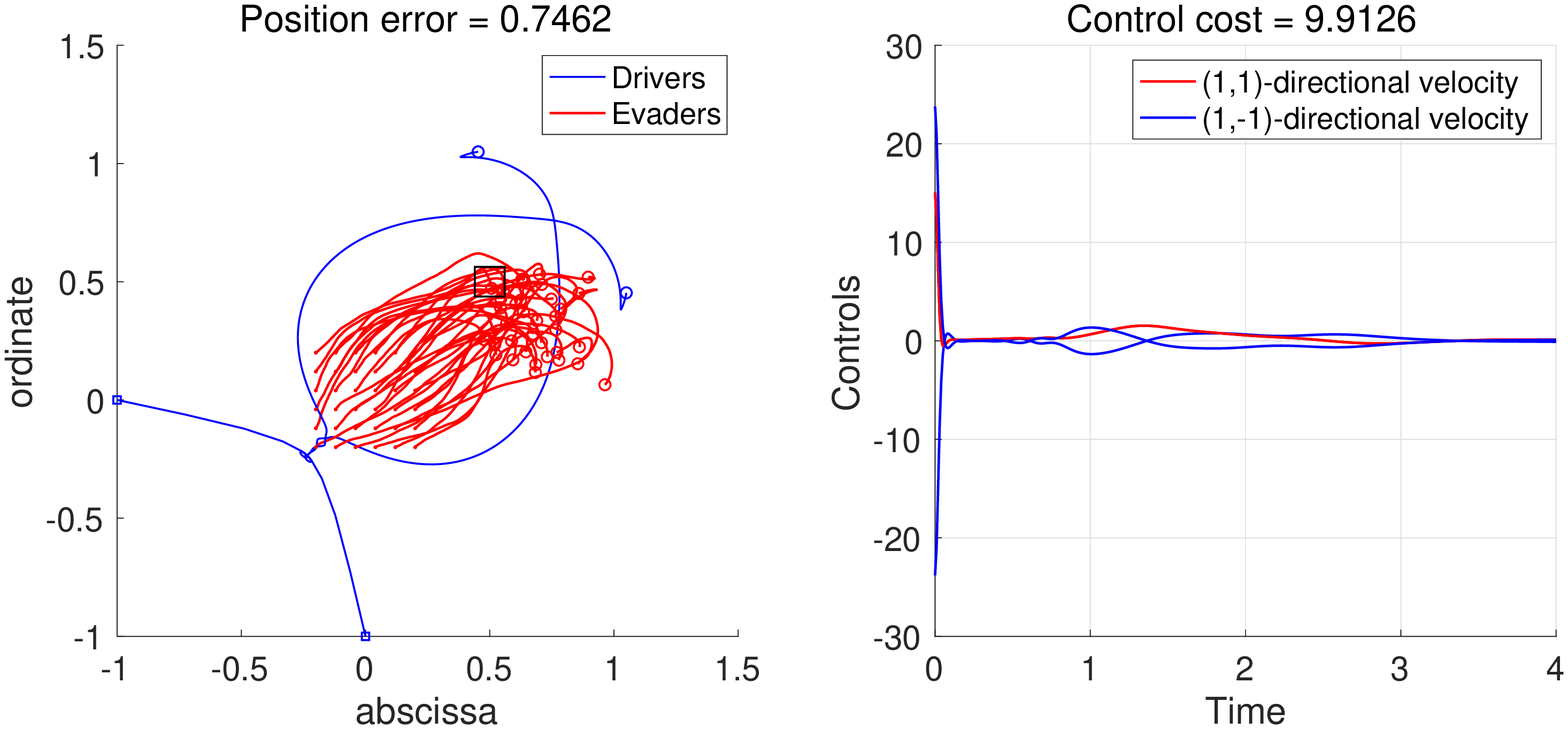}
    \includegraphics[width=0.8\textwidth]{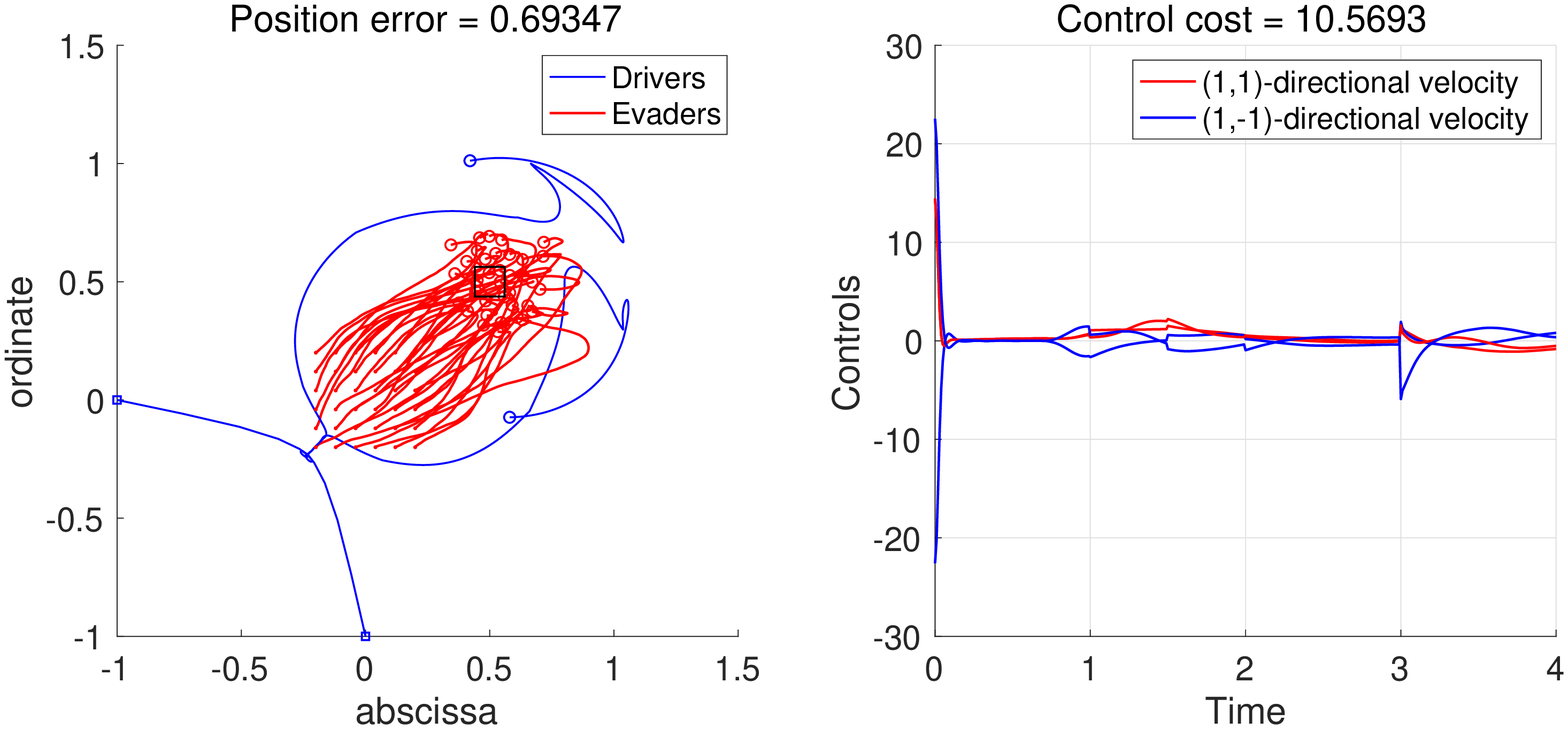}
  }
  \caption{The simulations on a noisy system with the open-loop control (above) and the control from the MPC-RBM algorithm with $P=2$, $\tau = 1$ and $\hat T = 3$ (below). The multiplicative noise is added on velocities of evaders, where the evaders escape to the right direction in the case of the open-loop control. \emph{Left}: the trajectories of the drivers and evaders for $t \in [0,4]$. \emph{Right}: the control functions along time.}
  \label{fig:test_noise}
\end{figure*}

\begin{table*}[t]
\begin{center}
\begin{tabular}{l|*{3}{c}|c}
\hline
& Control interval & GD iterations & Computation time & Cost $J$ \\
\hline
MPC-RBM, noisy system & $[0.0,1.0]$ & 1727 & \phantom{0}47.42 & 0.7040 \\
(Fig. \ref{fig:test_noise}) & $[1.0,2.0]$ & \phantom{0}429 & \phantom{0}12.27 & \\
(Total time: $124.33$)               & $[2.0,3.0]$ & \phantom{00}48 & \phantom{00}1.59 & \\
               & $[3.0,4.0]$ & 2252 & \phantom{0}63.05 & \\
\hline
\end{tabular}
\end{center}
\caption{The number of iterations and computation time (in seconds) for Fig. \ref{fig:test_noise} to find the optimal controls with the MPC-RBM algorithm.}\label{tab:time_mpc_noise}
\end{table*}

As we presented in Section \ref{sec:2-4}, the strategy of MPC is adopted to overcome the errors from the reduced dynamics. This effect can be seen significantly when the system has a noisy behavior. We consider a stochastic system by adding multiplicative noise ${\mathbf v}_i dB^i_t$ to the time derivatives of velocities in \eqref{eq:GBR}:
\begin{equation*}
\begin{aligned}
\begin{cases}
\dot{\mathbf x}_i = {\mathbf v}_i,\qquad i=1,\ldots,N,\\
d{\mathbf v}_i = \displaystyle\left[\frac{1}{N-1}\sum_{k=1,k\neq i}^N a({\mathbf x}_k-{\mathbf x}_i)({\mathbf v}_k-{\mathbf v}_i)
 - \frac{1}{N-1}\sum_{k=1,k\neq i}^N g({\mathbf x}_k-{\mathbf x}_i)({\mathbf x}_k-{\mathbf x}_i)\right.\\
\hspace{3.5em}\left.
- \displaystyle\frac{1}{M}\sum_{j=1}^M f({\mathbf y}_j-{\mathbf x}_i)({\mathbf y}_j-{\mathbf x}_i)\right]dt + \sigma {\mathbf v}_i dB^i_t, \qquad i = 1,\ldots,N\\
\dot {\mathbf y}_j = {\mathbf u}_j(t), \qquad j=1,\ldots,M,\\
{\mathbf x}_i(0) = {\mathbf x}_i^0,\quad {\mathbf v}_i(0) = {\mathbf v}_i^0, \quad {\mathbf y}_j(0) = {\mathbf y}_j^0,\qquad i=1,\ldots,N,\quad j=1,\ldots,M,
\end{cases}
\end{aligned}
\end{equation*}
where $(B^1_t,\ldots,B^N_t)$ is the $N$-dimensional Brownian motion and the noise strength $\sigma$ is a constant, $\sigma = 0.5$. The multiplicative noise is simulated with Milstein method.

Fig. \ref{fig:test_noise} shows the numerical simulations with the open-loop control (above) and the control from the MPC-RBM algorithm (below). The open-loop control is the one we calculated in Fig. \ref{fig:test1} from the original deterministic system \eqref{eq:GBR}. The MPC-RBM algorithm is operated from the RBM reduced model \eqref{eq:sGBR} with $\tau = 1$ and $\hat T = 3$ by updating the controlled trajectories of the noisy system. We set a smaller $\tau$ compared to the Fig. \ref{fig:test2} since the difference from the controlled system and the reduced model is increased. 

In the open-loop control, the evaders escape to an unexpected direction near the final time. Using the MPC-RBM algorithm, the drivers effectively surround the evaders with the information at $t=1$, $2$, and $3$. The computation time and the cost functions are described in Table \ref{tab:time_mpc_noise}, where the cost function is reduced to $0.7040$ from the cost of the open-loop control, $0.7561$.

\section{Conclusion and final remarks}\label{sec:final}

In this paper, we combine the model predictive control (MPC) with the random batch methods (RBM) to get a reliable control strategy in a short computation time. The suggested algorithm finds the optimal control on a reduced model from the RBM, as a predictive model in the process of MPC. 

The RBM simplifies the full system through the random sampling of the interactions. However, after the construction, the reduced model is a deterministic dynamics with a switching network in time. Hence, we may use various optimization algorithms to find the optimal control. 

With a given $P>1$, the computation cost on the dynamics is reduced to the order of $O(N(P+M))$ from $O(N(N+M))$ for the problem of $N$ evaders and $M$ drivers. Hence, in the MPC-RBM algorithm, the overall computation time is significantly reduced when there are plenty of evaders. 

One may choose a bigger $P$ to balance the approximation error with the computational cost. In the simulations of Fig. \ref{fig:error2}, increasing $P$ from $2$ to $4$ takes about $73\%$ more computations but has nearly $62\%$ smaller error in time evolution of positions. 

Though our focus is on a specific situation, the MPC-RBM algorithm works on a general interacting particle system. For example, the guiding problem can be simulated in a three-dimensional space or with a restriction on control. Still, the performance of the resulting control is not estimated rigorously due to the absence of analysis on the RBM and MPC in general complex systems.

\bibliographystyle{plain}
\bibliography{biblio}

\begin{thebibliography}{10}

\bibitem{Andersson2018}
J.~A.~E. Andersson, J.~Gillis, Horn G., J.~B. Rawlings, and M.~Diehl.
\newblock {CasADi} -- {A} software framework for nonlinear optimization and
  optimal control.
\newblock {\em Mathematical Programming Computation}, In Press, 2018.

\bibitem{bailo_optimal_2018}
R.~Bailo, M.~Bongini, J.~A. Carrillo, and D.~Kalise.
\newblock Optimal consensus control of the {Cucker}-{Smale} model.
\newblock {\em IFAC-PapersOnLine}, 51(13):1--6, 2018.

\bibitem{bongini_sparse_2015}
M.~Bongini, M.~Fornasier, O.~Junge, and B.~Scharf.
\newblock Sparse control of alignment models in high dimension.
\newblock {\em Networks and Heterogeneous Media}, 10(3):647--697, 2015.

\bibitem{bongini_mean-field_2017}
M.~Bongini, M.~Fornasier, F.~Rossi, and F.~Solombrino.
\newblock Mean-{Field} {Pontryagin} {Maximum} {Principle}.
\newblock {\em Journal of Optimization Theory and Applications}, 175(1):1--38,
  2017.

\bibitem{burger_controlling_2016}
M.~Burger, R.~Pinnau, A.~Roth, C.~Totzeck, and O.~Tse.
\newblock Controlling a self-organizing system of individuals guided by a few
  external agents -- particle description and mean-field limit.
\newblock {\em arXiv:1610.01325 [math]}, 2016.
\newblock arXiv: 1610.01325.

\bibitem{caponigro_sparse_2013}
M.~Caponigro, M.~Fornasier, B.~Piccoli, and E.~Trélat.
\newblock Sparse stabilization and optimal control of the {Cucker}-{Smale}
  model.
\newblock {\em Mathematical Control and Related Fields}, 3(4):447--466, 2013.

\bibitem{carrillo_sharp_2017}
J.~A. Carrillo, Y.-P. Choi, P.~B. Mucha, and J.~Peszek.
\newblock Sharp conditions to avoid collisions in singular {Cucker}–{Smale}
  interactions.
\newblock {\em Nonlinear Analysis: Real World Applications}, 37:317--328, 2017.

\bibitem{carrillo_analytical_2018}
J.~A. Carrillo, Y.-P. Choi, C.~Totzeck, and O.~Tse.
\newblock An analytical framework for consensus-based global optimization
  method.
\newblock {\em Mathematical Models and Methods in Applied Sciences},
  28(06):1037--1066, 2018.

\bibitem{carrillo_consensus-based_2019}
J.~A. Carrillo, S.~Jin, L.~Li, and Y.~Zhu.
\newblock A consensus-based global optimization method for high dimensional
  machine learning problems.
\newblock {\em arXiv:1909.09249 [math]}, 2019.
\newblock arXiv: 1909.09249.

\bibitem{cucker_avoiding_2010}
F.~Cucker and J.-G. Dong.
\newblock Avoiding {Collisions} in {Flocks}.
\newblock {\em IEEE Transactions on Automatic Control}, 55(5):1238--1243, 2010.

\bibitem{cucker_emergent_2007}
F.~Cucker and S.~Smale.
\newblock Emergent {Behavior} in {Flocks}.
\newblock {\em IEEE Transactions on Automatic Control}, 52(5):852--862, 2007.

\bibitem{dorfler_synchronization_2012}
F.~Dorfler and F.~Bullo.
\newblock Synchronization and transient stability in power networks and
  nonuniform kuramoto oscillators.
\newblock {\em SIAM Journal on Control and Optimization}, 50(3):1616--1642,
  2012.

\bibitem{escobedo_optimal_2016}
R.~Escobedo, A.~Ibañez, and E.~Zuazua.
\newblock Optimal strategies for driving a mobile agent in a “guidance by
  repulsion” model.
\newblock {\em Communications in Nonlinear Science and Numerical Simulation},
  39:58--72, 2016.

\bibitem{fax2004information}
J~Alexander Fax and Richard~M Murray.
\newblock Information flow and cooperative control of vehicle formations.
\newblock {\em IEEE transactions on automatic control}, 49(9):1465--1476, 2004.

\bibitem{gade_herding_2015}
S.~Gade, A.~A. Paranjape, and S.-J. Chung.
\newblock Herding a {Flock} of {Birds} {Approaching} an {Airport} {Using} an
  {Unmanned} {Aerial} {Vehicle}.
\newblock In {\em {AIAA} {Guidance}, {Navigation}, and {Control} {Conference}},
  Kissimmee, Florida, 2015. American Institute of Aeronautics and Astronautics.

\bibitem{garcia_model_1989}
C.~E. García, D.~M. Prett, and M.~Morari.
\newblock Model predictive control: {Theory} and practice—{A} survey.
\newblock {\em Automatica}, 25(3):335--348, 1989.

\bibitem{golse_mean-field_2003}
F.~Golse.
\newblock The mean-field limit for the dynamics of large particle systems.
\newblock {\em Journées équations aux dérivées partielles}, pages 1--47,
  2003.

\bibitem{grune2017nonlinear}
L.~Gr{\"u}ne and J.~Pannek.
\newblock Nonlinear model predictive control.
\newblock In {\em Nonlinear Model Predictive Control}, pages 45--69. Springer,
  2017.

\bibitem{ha_convergence_2019}
S.-Y. Ha, S.~Jin, and D.~Kim.
\newblock Convergence of a first-order consensus-based global optimization
  algorithm.
\newblock {\em arXiv:1910.08239 [math]}, 2019.

\bibitem{ha2009simple}
S.-Y. Ha and J.-G. Liu.
\newblock A simple proof of the cucker-smale flocking dynamics and mean-field
  limit.
\newblock {\em Communications in Mathematical Sciences}, 7(2):297--325, 2009.

\bibitem{jin_random_2020}
S.~Jin, L.~Li, and J.-G. Liu.
\newblock Random {Batch} {Methods} ({RBM}) for interacting particle systems.
\newblock {\em Journal of Computational Physics}, 400:108877, 2020.

\bibitem{ko2019GBR}
D.~Ko and E.~Zuazua.
\newblock Asymptotic behavior and control of a ``guidance by repulsion" model.
\newblock {\em Mathematical Models and Methods in Applied Sciences}, Online
  Ready. doi: 10.1142/S0218202520400047, 2020.

\bibitem{jyh-ming_lien_shepherding_2004}
J.~M. Lien, O.~B. Bayazit, R.~T. Sowell, S.~Rodriguez, and N.~M. Amato.
\newblock Shepherding behaviors.
\newblock In {\em {IEEE} {International} {Conference} on {Robotics} and
  {Automation}, 2004. {Proceedings}. {ICRA} '04. 2004}, pages 4159--4164 Vol.4,
  New Orleans, LA, USA, 2004. IEEE.

\bibitem{ma_finite-time_2017}
J.~Ma and E.~M.-K. Lai.
\newblock Finite-time flocking control of a swarm of cucker-smale agents with
  collision avoidance.
\newblock In {\em 2017 24th {International} {Conference} on {Mechatronics} and
  {Machine} {Vision} in {Practice} ({M2VIP})}, pages 1--6, Auckland, 2017.
  IEEE.

\bibitem{mhaskar2006robust}
P.~Mhaskar.
\newblock Robust model predictive control design for fault-tolerant control of
  process systems.
\newblock {\em Industrial \& engineering chemistry research},
  45(25):8565--8574, 2006.

\bibitem{motsch_new_2011}
S.~Motsch and E.~Tadmor.
\newblock A {New} {Model} for {Self}-organized {Dynamics} and {Its} {Flocking}
  {Behavior}.
\newblock {\em Journal of Statistical Physics}, 144(5):923--947, 2011.

\bibitem{nesterov2013introductory}
Y.~Nesterov.
\newblock {\em Introductory lectures on convex optimization: A basic course},
  volume~87.
\newblock Springer Science \& Business Media, 2013.

\bibitem{nikolaou2001model}
M.~Nikolaou.
\newblock Model predictive controllers: A critical synthesis of theory and
  industrial needs.
\newblock {\em Advances in Chemical Engineering}, 26:131--204, 2001.

\bibitem{park_cucker-smale_2010}
J.~Park, H.~J. Kim, and S.-Y. Ha.
\newblock Cucker-{Smale} {Flocking} {With} {Inter}-{Particle} {Bonding}
  {Forces}.
\newblock {\em IEEE Transactions on Automatic Control}, 55(11):2617--2623,
  2010.

\bibitem{piccoli_sparse_2019}
B.~Piccoli, N.~P. Duteil, and E.~Trélat.
\newblock Sparse {Control} of {Hegselmann}--{Krause} {Models}: {Black} {Hole}
  and {Declustering}.
\newblock {\em SIAM Journal on Control and Optimization}, 57(4):2628--2659,
  2019.

\bibitem{pinnau_interacting_2018}
R.~Pinnau and C.~Totzeck.
\newblock Interacting {Particles} and {Optimization}.
\newblock {\em PAMM}, 18(1), 2018.

\bibitem{pontryagin2018mathematical}
L.~S. Pontryagin.
\newblock {\em Mathematical theory of optimal processes}.
\newblock Routledge, 2018.

\bibitem{porfiri2008criteria}
M.~Porfiri and M.~Di~Bernardo.
\newblock Criteria for global pinning-controllability of complex networks.
\newblock {\em Automatica}, 44(12):3100--3106, 2008.

\bibitem{prett1987design}
D.~M. Prett and C.~E. Garcia.
\newblock Design of robust process controllers.
\newblock {\em IFAC Proceedings Volumes}, 20(5):275--280, 1987.

\bibitem{strombom_solving_2014}
D.~Strömbom, R.~P. Mann, A.~M. Wilson, S.~Hailes, A.~J. Morton, D.~J.~T.
  Sumpter, and A.~J. King.
\newblock Solving the shepherding problem: heuristics for herding autonomous,
  interacting agents.
\newblock {\em Journal of The Royal Society Interface}, 11(100):20140719, 2014.

\bibitem{tanner_flocking_2007}
H.~G. Tanner, A.~Jadbabaie, and G.~J. Pappas.
\newblock Flocking in {Fixed} and {Switching} {Networks}.
\newblock {\em IEEE Transactions on Automatic Control}, 52(5):863--868, 2007.

\bibitem{trelat2005controle}
E.~Tr{\'e}lat.
\newblock {\em Contr{\^o}le optimal: th{\'e}orie \& applications}.
\newblock Vuibert Paris, 2005.

\end{thebibliography}

\end{document}